\documentclass[12pt,a4paper]{amsart}
\usepackage{amssymb}
\usepackage[a4paper,centering,twoside,textwidth=6in]{geometry}
\usepackage[arrow,matrix]{xy}

\theoremstyle{plain}
\newtheorem{theorem}{Theorem}[section]
\newtheorem{conj}[theorem]{Conjecture}
\newtheorem{corollary}[theorem]{Corollary}
\newtheorem{prop}[theorem]{Proposition}
\newtheorem{Lemma}[theorem]{Lemma}
\newtheorem{claim}[theorem]{Claim}

\theoremstyle{definition}
\newtheorem{defn}[theorem]{Definition}

\newtheorem{remark}[theorem]{Remark}

\newcommand{\thmref}[1]{Theorem~\ref{#1}}
\newcommand{\secref}[1]{Section~\ref{#1}}

\newcommand{\lemref}[1]{Lemma~\ref{#1}}
\newcommand{\corref}[1]{Corollary~\ref{#1}}

\newcommand{\remarkref}[1]{Remark~\ref{#1}}

\newcommand{\propref}[1]{Proposition~\ref{#1}}

\DeclareMathOperator{\PSL}{PSL} \DeclareMathOperator{\SL}{SL}
\DeclareMathOperator{\PGL}{PGL} \DeclareMathOperator{\GL}{GL}
\DeclareMathOperator{\tr}{tr} \DeclareMathOperator{\lcm}{lcm}
\newcommand{\BZ}{\mathbb{Z}}

\newcommand{\BP}{\mathbb{P}}
\newcommand{\BF}{\mathbb{F}}
\newcommand{\BA}{\mathbb{A}}
\newcommand{\ep}{\epsilon}
\newcommand{\al}{\alpha}
\newcommand{\be}{\beta}
\newcommand{\ga}{\gamma}
\newcommand{\de}{\delta}
\newcommand{\la}{\lambda}
\newcommand{\om}{\omega}
\newcommand{\vp}{\varphi}

\newcommand{\fc}{\frac}

\begin{document}

\title {Surjectivity and equidistribution of the word $x^ay^b$ on $\PSL(2,q)$ and $\SL(2,q)$}
\author{Tatiana Bandman and Shelly Garion}

\address{Bandman: Department of Mathematics, Bar-Ilan University, 52900 Ramat Gan, ISRAEL}
\email{bandman@macs.biu.ac.il}

\address{Garion: Institut des Hautes \'{E}tudes Scientifiques, 91440 Bures-sur-Yvette, FRANCE}
\email{shellyg@ihes.fr}



\maketitle

\begin{abstract}
We determine the integers $a,b \geq 1$ and the prime powers $q$ for
which the word map $w(x,y)=x^ay^b$ is surjective on the group
$\PSL(2,q)$ (and $\SL(2,q)$). We moreover show that this map is almost
equidistributed for the family of groups $\PSL(2,q)$ (and $\SL(2,q)$).
Our proof is based on the investigation of the trace map of positive words.
\end{abstract}

\section{Introduction}

\subsection{Word maps in finite simple groups}
During the last years there has been a great interest in \emph{word
maps} in groups, for an extensive survey see \cite{Se}. These maps
are defined as follows. Let $w=w(x_1,\dots,x_d)$ be a non-trivial
\emph{word}, namely a non-identity element of the free group $F_d$
on the generators $x_1,\dots,x_d$. Then we may write
$w=x_{i_1}^{n_1}x_{i_2}^{n_2}\dots x_{i_k}^{n_k}$ where $1 \leq i_j
\leq d$, $n_j \in \BZ$, and we may assume further that $w$ is
reduced. Let $G$ be a group. For $g_1,\dots,g_d$ we write
$$
   w(g_1,\dots,g_d)=g_{i_1}^{n_1}g_{i_2}^{n_2}\dots g_{i_k}^{n_k} \in G,
$$
and define
$$
   w(G) = \{w(g_1,\dots,g_d): g_1,\dots,g_d \in G \},
$$
as the set of values of $w$ in $G$. The corresponding map $w:G^d
\rightarrow G$ is called a \emph{word map}.

Borel \cite{Bo} showed that the word map induced by $w \neq 1$ on
simple algebraic groups is a dominant map. Larsen \cite{La} used
this result to show that for every non-trivial word $w$ and $\ep>0$
there exists a number $C(w,\ep)$ such that if $G$ is a finite simple
group with $|G|>C(w,\ep)$ then $|w(G)| \geq |G|^{1-\ep}$. By recent
work of Larsen, Shalev and Tiep \cite{Sh09,LS,LST}, for every non-trivial
word $w$ there exists a constant $C(w)$ such that if $G$ is a finite
simple group satisfying $|G|>C(w)$ then $w(G)^2=G$.

It is therefore interesting to find words $w$ for which $w(G)=G$ for
any finite simple non-abelian group $G$. Due to immense work spread
over more than 50 years, it is now known that the commutator word
$w=[x,y] \in F_2$ satisfies $w(G)=G$ for any finite simple group $G$
(see \cite{LOST10} and the references therein). On the other hand, it
is easy to see that if $G$ is a finite group and $b$ is an integer
which is not relatively prime to the order of $G$ then $w(G) \neq G$
for the word $w=x^b$.

The words of the form $w=x^ay^b \in F_2$ have also attracted special
interest. Due to recent work of Larsen, Shalev and Tiep \cite{LST},
it is known that any such word is surjective on sufficiently large
finite simple groups (see \cite[Theorem 1.1.1 and Corollary
1.1.3]{LST}), more precisely,

\begin{theorem}\cite{LST}.\label{thm.LST}
Let $a,b$ be two non-zero integers. Then there exists a number $N =
N(a,b)$ such that if $G$ is a finite simple non-abelian group of
order at least $N$, then any element in $G$ can be written as
$x^ay^b$ for some $x,y \in G$.
\end{theorem}

By further recent results of Guralnick and Malle \cite{GM}
and of Liebeck, O`Brien, Shalev and Tiep \cite{LOST11},
some words of the form $x^by^b$ are known to be surjective on \emph{all}
finite simple groups.

\begin{theorem}\cite[Corollary 1.5]{GM}. \label{thm.GM}
Let $G$ be a finite simple non-abelian group and let $b$ be either a
prime power or a power of $6$. Then any element in $G$ can be
written as $x^by^b$ for some $x,y \in G$.
\end{theorem}

Note that in general, the word $x^by^b$ is not necessarily
surjective on \emph{all} finite simple groups. Indeed, if $b$ is a
multiple of the exponent of $G$ then necessarily $x^by^b=id$ for
every $x,y \in G$. It is therefore interesting to find more examples
for words of the form $x^by^b$ which are \emph{not} surjective on
all finite simple groups.

\medskip

More generally, one can ask whether it is possible to generalize
\thmref{thm.LST} for other word maps. In particular, the following conjecture was
raised:

\begin{conj}[Shalev]\cite[Conjectures 2.8 and 2.9]{Sh07}. \label{conj.surjall}
Let $w \ne 1$ be a word which is not a proper power of another word.
Then there exists a number $C(w)$ such that, if $G$ is either $A_r$
or a finite simple group of Lie type of rank $r$, where $r > C(w)$,
then $w(G)=G$.
\end{conj}

It is also interesting to investigate the \emph{distribution} of the word
map. Due to the work of Garion and Shalev \cite{GS}, it is known
that the word $w=x^2y^2$ is \emph{almost equidistributed} for the
family of finite simple groups, namely,

\begin{theorem}\cite[Theorem 7.1]{GS}.
Let $G$ be a finite simple group, and let $ w:G \times G \to G$ be
the map given by $w(x,y) = x^2y^2$. Then there is a subset $S
\subseteq G$ with $|S| = (1-o(1))|G|$ such that $|w^{-1}(g)| =
(1+o(1))|G|$ for all $g \in S$.
Where $o(1)$ denotes a function depending only on $G$ which tends to zero as $|G| \rightarrow \infty$.
\end{theorem}

Another question which was raised by Shalev \cite[Problem
2.10]{Sh07} is which words $w$ induce an \emph{almost
equidistributed} map for the family of finite simple groups. In
particular, does words of the form $w=x^ay^b$ induce \emph{almost
equidistributed} maps?

\subsection{The word $w(x,y)=x^ay^b$ on the groups $\SL(2,q)$ and $\PSL(2,q)$}
In this paper we analyze the word map $x^ay^b$ in the groups
$\SL(2,q)$ and $\PSL(2,q)$. Analysis of Engel word maps in these
groups was carried out in our previous work \cite{BGG}.

We start by determining precisely the positive integers $a,b$ and
prime powers $q$ for which the word map $w=x^ay^b$ is surjective on
$\SL(2,q) \setminus \{-id\}$ and on $\PSL(2,q)$.

\begin{defn}\label{def.non.degenerate}
Let $a,b \geq 1$ and let $q=p^e$ be a prime power. We say that the
word $w=x^ay^b$ is \emph{non-degenerate} with respect to $q$ if and
only if none of the following conditions holds:
\begin{itemize}
\item $p=2$, $a$ is a multiple of $2(q^2-1)$ and $b$ is not relatively prime to $2(q^2-1)$;
\item $p=2$, $b$ is a multiple of $2(q^2-1)$ and $a$ is not relatively prime to $2(q^2-1)$;
\item $p$ is odd, $a$ is a multiple of $\frac{p(q^2-1)}{4}$ and $b$ is not
relatively prime to $\frac{p(q^2-1)}{4}$;
\item $p$ is odd, $b$ is a multiple of $\frac{p(q^2-1)}{4}$ and $a$ is not
relatively prime to $\frac{p(q^2-1)}{4}$.
\end{itemize}
\end{defn}

Obviously, if the word map $w=x^ay^b$ is surjective on $\PSL(2,q)$
then it is necessarily non-degenerate with respect to $q$. On the
other hand, we prove the following proposition.

\begin{prop}\label{prop.abne2}
If $w=x^ay^b$ is non-degenerate with respect to $q$, then
all \emph{semisimple} elements, namely matrices in $\SL(2,q)$ whose
trace is not $\pm 2$, are in the image of the word map $w$.
\end{prop}

Unfortunately, even if $w=x^ay^b$ is non-degenerate with respect to
$q$, the image of $w=x^ay^b$ may not contain the \emph{unipotent}
elements, namely, matrices $\pm id \ne z \in \SL(2,q)$ satisfying
$\tr(z) = \pm 2$. This phenomenon happens when one of the following
obstructions occurs.

\begin{defn}\label{def.obs}
Let $a,b \geq 1$ and let $q$ be a prime power.
We define the following \emph{obstructions}:
\begin{itemize}
\item \emph{Obstruction (i):} $q=2^e$, $e$ is odd, and $a,b$ are divisible by $\frac{2(q^2-1)}{3}$;
\item \emph{Obstruction (ii):} $q \equiv 3 \bmod 4$ and $a,b$ are divisible by $\frac{p(q^2-1)}{8}$;
\item \emph{Obstruction (iii):} $q \equiv 11 \bmod 12$ and $a,b$ are divisible by $\frac{p(q^2-1)}{6}$;
\item \emph{Obstruction (iv):} $q \equiv 5 \bmod 12$ and $a,b$ are divisible by $\frac{p(q^2-1)}{12}$.
\end{itemize}
\end{defn}

\begin{theorem}\label{thm.surj.SL.even}
Let $e \geq 1$ and let $q=2^e$. Let $a,b \geq 1$.

Then the word map $w=x^ay^b$ is surjective on $\SL(2,q) = \PSL(2,q)$
if and only if $w$ is non-degenerate with respect to $q$ and
obstruction (i) does not occur.
\end{theorem}

\begin{theorem}\label{thm.surj.SL.odd}
Let $p$ be an odd prime number, $e \geq 1$ and $q=p^e$.
Let $a,b \geq 1$.

Then the word map $w=x^ay^b$ is surjective on $\SL(2,q) \setminus
\{-id\}$ if and only if $w$ is non-degenerate with respect to $q$
and none of the obstructions (ii),(iii),(iv) occurs.
\end{theorem}

\begin{theorem}\label{thm.surj.PSL}
Let $p$ be an odd prime number, $e \geq 1$ and $q=p^e$.
Let $a,b \geq 1$.

Then the word map $w=x^ay^b$ is surjective on $\PSL(2,q)$ if and
only if $w$ is non-degenerate with respect to $q$ and obstruction (ii)
does not occur.
\end{theorem}

For example, we deduce that the word $w=x^{42}y^{42}$ is \emph{not}
surjective on the groups $\PSL(2,7)$ and $\PSL(2,8)$.

The last theorem implies that for the family of groups $\PSL(2,q)$
one can give a precise estimation for the bound $N=N(a,b)$ appearing
in \thmref{thm.LST}.
\begin{corollary}\label{cor.Nab}
For every $a,b \geq 1$, let
$$ Q=Q(a,b) = \max \{\sqrt{3a}, \sqrt{3b}\}, $$
$$ N=N(a,b) = \max \bigl\{\frac{3\sqrt{3}}{2}a^{3/2}, \frac{3\sqrt{3}}{2}b^{3/2}\bigr\} .$$
Then the word map $w=x^ay^b$ is surjective on the group $\PSL(2,q)$ for any
$q>Q$, and hence whenever $|\PSL(2,q)|>N$.
\end{corollary}

However, the statement of \thmref{thm.LST} \cite{LST} no longer
holds for the quasi-simple group $\SL(2,q)$, as indicated by the
following theorem and its corollary.

\begin{theorem}\label{thm.min.id}
Let $q$ be an odd prime power, and set $K = \max \left\{k: 2^k \
divides \ \frac{q^2-1}{2} \right\}.$ Let $a,b \geq 1$. Then $-id
\neq x^ay^b$ for every $x,y \in \SL(2,q)$ if and only if $2^K$
divides both $a$ and $b$.
\end{theorem}

\begin{corollary}\label{cor.4.4}
If $q \equiv \pm 3 \bmod 8$, then $x^4y^4 \ne -id$ for every $x,y
\in \SL(2,q)$.
\end{corollary}

In addition, we show that for any $a,b \geq 1$, the word map
$w=x^ay^b$ is \emph{almost equiditributed} for the family of groups
$\PSL(2,q)$ (and $\SL(2,q)$).

\begin{theorem}\label{thm.equi}
Let $q$ be a prime power and let $G$ be either the group
$\SL(2,q)$ or the group $\PSL(2,q)$.

Let $a,b \geq 1$ and let $ w:G \times G \to G$ be
the map given by $w(x,y) = x^ay^b$. Then there is a subset $S
\subseteq G$ with $|S| = (1-o(1))|G|$ such that $|w^{-1}(g)| =
(1+o(1))|G|$ for all $g \in S$.
Where $o(1)$ denotes a function of $q$ which tends to zero as $q \rightarrow \infty$.
\end{theorem}

\subsection{Organization and outline of the proof}
For the convenience of the reader, we describe the organization of
the paper, as well as give a bird's eye view of the proofs.

In \secref{sect.trace.map} we compute the \emph{trace map} of the
word $w(x,y)=x^ay^b$ (\lemref{lem.ab}), and more generally, of any
\emph{positive} word in $F_2$ (\thmref{thm.gen.word}). For any word
$w=w(x,y) \in F_2$, the trace map $\tr(w)$ is a polynomial
$P(s,u,t)$ in $s=\tr(x), t=\tr(y)$ and $u=\tr(xy)$.

In \secref{sect.grp} we collect basic facts on the surjectivity of
$w=x^ay^b$ on finite groups in general, and in \secref{sect.prop.SL}
we describe some properties of the groups $\SL(2,q)$ and $\PSL(2,q)$
that are used later on.

By \lemref{lem.ab}, $\tr(x^ay^b)$ is a \emph{linear} polynomial in
$u$. We deduce in \secref{sect.surj.trace.map} that if neither $a$
nor $b$ is divisible by the exponent of $\PSL(2,q)$, then any
element in $\BF_q$ can be written as $\tr(x^ay^b)$ for some $x,y \in
\SL(2,q)$ (\corref{cor.tr.surj}). This immediately implies
\propref{prop.abne2}, stating that if $w=x^ay^b$ is non-degenerate
with respect to $q$, any \emph{semisimple} element (namely, $z \in
\SL(2,q)$ with $\tr(z) \neq \pm 2$) can be written as $z=x^ay^b$ for
some $x,y \in \SL(2,q)$.

However, when $z$ is \emph{unipotent} (namely, $z \neq \pm id$ and
$\tr(z) = \pm 2$) one has to be more careful, and a detailed
analysis is done in \secref{sect.non.surj}. Indeed, it may happen
that $w=x^ay^b$ is non-degenerate with respect to $q$, but
nevertheless the image of the word map $w=x^ay^b$ does not contain
any unipotent (see Propositions \ref{prop.ab2} and
\ref{prop.abmin2}).

These are the ingredients needed for the proofs of Theorems
\ref{thm.surj.SL.even}, \ref{thm.surj.SL.odd} and
\ref{thm.surj.PSL}, on the surjectivity of the word $w=x^ay^b$ on
$\PSL(2,q)$ and $\SL(2,q) \setminus \{-id\}$, which are presented in
\secref{sect.surj}. In addition, we determine in
\secref{sect.min.id} when $-id$ can be written as $x^ay^b$ for some
$x,y \in \SL(2,q)$, thus proving \thmref{thm.min.id}.

In \secref{sect.equi} we prove \thmref{thm.equi} and show that the
word map $w=x^ay^b$ is almost equidistributed for the family of
groups $\PSL(2,q)$ (and $\SL(2,q)$). The basic idea is
to show that for a general $\al \in \BF_q$, the surface
$$S_\al(\BF_q)=\{P(s,u,t)=\tr(z)=\al\}\subset \BA^3(\BF_q)$$
is birational to a plane $\BA^2_{s,t}$. As a result, we are able not only to find points on
$S_\al(\BF_q),$  but even to estimate their number.

\begin{remark}
After this paper was completed, we were informed that M. Larsen and
A. Shalev have proved the uniformity of such words in more general
contexts. In another work, M. Larsen, A. Shalev and P.H. Tiep
\cite{LST11} have showed non-surjectivity of some words in
quasi-simple groups, and in particular obtained the same result of
\corref{cor.4.4}.
\end{remark}

\subsection*{Acknowledgement}
Bandman is supported in part by Ministry of Absorption (Israel),
Israeli Academy of Sciences and Minerva Foundation (through the Emmy
Noether Research Institute of Mathematics).

Garion is supported by a European Post-doctoral Fellowship (EPDI),
during her stay at the Institut des Hautes \'{E}tudes Scientifiques
(Bures-sur-Yvette).

The authors are grateful to A. Shalev for discussing his questions and
conjectures with them. They are also thankful for B. Kunyavskii, M. Larsen and
A. Mann for useful discussions.

\section{The trace map}\label{sect.trace.map}
\subsection{The trace map}
The \emph{trace map} method is based on the following classical
Theorem (see, for example, \cite{Vo,Fr,FK} or \cite{Mag, Go} for a
more modern exposition).

\begin{theorem}[Trace map]\label{thm.trace.map}
Let $F=\left<x,y\right>$ denote the free group on two generators.
Let us embed $F$ into $\SL(2,\BZ)$ and denote by $\tr$ the trace
character. If $w$ is an arbitrary element of $F$, then the character
of $w$ can be expressed as a polynomial
$$\tr(w)=P(s,u,t)$$
with integer coefficients in the three characters $s=\tr(x),
u=\tr(xy)$ and $t=\tr(y)$.
\end{theorem}

Note that the same remains true for the group $\SL(2,q)$. The
general case, $\SL(2,R)$, where $R$ is a commutative ring, can be
found in \cite{CMS}.

The following theorem is originally due to Macbeath~\cite{Mac} and
was used by Bandman, Grunewald, Kunyavskii and Jones to investigate
verbal dynamical systems in the group $\SL(2,q)$ (see ~\cite[Theorem
3.4]{BGK}).

\begin{theorem}\cite[Theorem 1]{Mac}. \label{thm.Mac}
For any $(s,u,t) \in \BF_q^3$ there exist two matrices $x,y \in
\SL(2,q)$ satisfying $\tr(x)=s, \tr(y)=t$ and $\tr(xy)=u$.
\end{theorem}

\subsection{Trace map of the word $w(x,y)=x^ay^b$}\label{sect.trace.xatb}

The following Lemma shows that the trace map of the word $w(x,y)=x^ay^b$ is a
linear polynomial in $\tr(xy)$.

\begin{Lemma}\label{lem.ab}
Let $w(x,y)=x^ay^b$ where $a,b\geq 1$ and $x,y \in \SL(2,q).$ Let
$s=\tr(x), u=\tr(xy), t=\tr(y).$ Then $$\tr  w(x,y)= u\cdot
f_{a,b}(s,t)+h_{a,b}(s,t),$$ where
$$f_{a,b}(s,t), \ h_{a,b}(s,t)\in \BF_q[s,t]$$ are polynomials satisfying:

\begin{itemize}
\item the highest degree summand of $f_{a,b}(s,t)$ (of degree $a+b-2$) is:
$$s^{a-1}t^{b-1};$$
\item the highest degree summand of $h_{a,b}(s,t)$ (of degree $a+b-2$) is:
\[
\begin{cases}
- (s^{a}t^{b-2} + s^{a-2}t^{b}) & \text{if } a>1 \text{ and } b>1; \\
- st^{b-2} & \text{if } a=1 \text{ and } b>1; \\
- s^{a-2}t & \text{if } a>1 \text{ and } b=1.
\end{cases}
\]
\end{itemize}
\end{Lemma}

\begin{proof}
We need to prove that the polynomials $f(s,t)=f_{a,b}(s,t)$ and
$h(s,t)=h_{a,b}(s,t)$ satisfy the following properties:

\begin{enumerate}\renewcommand{\theenumi}{\it \roman{enumi}}
\item the coefficient of $u,$ $f(s,t),$ has precisely one monomial
summand $s^{a-1}t^{b-1}$ with coefficient $1;$
\item for all other monomial summands $c_{i,j} s^it^j, \ c_{i,j}\in\BF_q, $
of $f(s,t),$ the following inequalities hold: $i\leq a-1,$  $j\leq
{b-1},$ and $i+j< a+b-2;$
\item $h(s,t)$ contains the summand
$(s^{a}t^{b-2} + s^{a-2}t^{b})$ with coefficient $-1$;
\item for all other monomial summands $c_{i,j} s^it^j, \ c_{i,j}\in\BF_q, $
of $h(s,t),$ the following inequalities hold: $i\leq a-2,$  $j\leq
{b-2},$ and so $i+j\leq a+b-4.$
\end{enumerate}

We prove these properties by induction on $a+b$, using the
well-known formula
\begin{equation}\label{tr.formula}
\tr(AB)+\tr(AB^{-1})=\tr(A)\tr(B).
\end{equation}

\emph{Induction base.} $a\leq 3,b\leq 3 .$ In these cases,
$$\tr(xy)=u, \
\tr(x^2y)=us-t, \ \tr(xy^2)=ut-s, \ \tr(x^2y^2)= ust-s^2-t^2+2,$$
$$\tr(xy^3) =  (t^2-1)u - st, \ \ \tr(x^3y) = (s^2-1)u - st,$$
$$\tr(x^2y^3) = (st^2 - s)u - s^2t - t^3 + 3t, \ \tr(x^3y^2) = (s^2t
- t)u - s^3 - st^2 + 3s,$$
$$\tr(x^3y^3)= (s^2t^2 - s^2 - t^2 + 1)u - s^3t - st^3 + 4st.$$

\emph{Induction hypothesis.} Assume that the Lemma is valid for
$a+b< n$ for some $n \geq 5$.

\emph{Induction Step.} We prove the claim for $a+b=n$, by
considering the following cases:

\medskip

{\bf  Case  1.}  $w(x,y)=x^ay^b,  \ a \geq b$ .

Using \eqref{tr.formula} we get:
\begin{align*}
\tr (x^ay^b) &= \tr(x)\tr (x^{a-1}y^b)- \tr (xy^{-b}x^{-a+1}) \\
&= \tr(x)\tr (x^{a-1}y^b)- \tr (x^{a-2}y^b) \\
&= s(u\cdot f_{a-1,b}(s,t)+h_{a-1,b}(s,t))-(u\cdot f_{a-2,b}(s,t)+h_{a-2,b}(s,t)) \\
&= u(s \cdot f_{a-1,b}(s,t)-f_{a-2,b}(s,t)) + (s\cdot h_{a-1,b}(s,t))-h_{a-2,b}(s,t)) \\
&= u\cdot f_{a,b}(s,t)+h_{a,b}(s,t).
\end{align*}

By the induction hypothesis, the resulting polynomial:
\begin{enumerate}\renewcommand{\theenumi}{\it \roman{enumi}}
\item is linear in $u;$

\item the highest degree summand of $f_{a,b}(s,t)$ is:
$$ss^{(a-1)-1}t^{b-1} = s^{a-1}t^{b-1} ;$$

\item for all other monomial summands $c_{i,j} s^it^j$ of $f(s,t)$ the following inequalities hold:
$i \leq (a-1)-1+1 = a-1$, $j \leq b-1$, and $$i+j < (a-1)+1+b-2=
a+b-2;$$

\item the highest degree summand in $h(s,t)$ is:
$$-s (s^{a-1}t^{b-2} + s^{(a-1)-2}t^b) = - (s^{a}t^{b-2} + s^{a-2}t^b).$$

\item for all other monomial summands $c_{i,j} s^it^j$ of $h(s,t)$ the
following inequalities hold: $i \leq (a-1)-2+1 = a-2$, $j \leq b-2$,
and so $$i+j \leq (a-1)+1+b-4= a+b-4;$$
\end{enumerate}

{\bf  Case  2.}  $w(x,y)=x^ay^b,  \ a<b$.

Similarly we get:
\begin{align*}
\tr (x^ay^b) &= \tr (x^ay^{b-1})\tr(y) - \tr(x^ay^{b-1}y^{-1}) \\
&= \tr(y)\tr (x^ay^{b-1}) - \tr(x^ay^{b-2}) \\
&= t(u\cdot f_{a,b-1}(s,t)+h_{a,b-1}(s,t)) - (u\cdot f_{a,b-2}(s,t)+h_{a,b-2}(s,t))= \\
&= u(t\cdot f_{a,b-1}(s,t) + f_{a,b-2}(s,t)) + (t\cdot h_{a,b-1}(s,t) - h_{a,b-2}(s,t))= \\
&= u\cdot f_{a,b}(s,t) + h_{a,b}(s,t).
\end{align*}

Similarly to Case 1 we get a polynomial satisfying the desired properties.
\end{proof}

\begin{remark}\label{rem.non-positive-1}
Assume that $a,b\ne0$ but not necessarily positive. Since
$\tr(xy^{-1})=st-u,$ we deduce from \lemref{lem.ab} that
$$\tr(x^ay^b)=u\cdot f_{a,b}(s,t)+h_{a,b}(s,t),$$ where
the highest degree summand of $f_{a,b}(s,t)$ (of degree $a+b-2$) is
$\pm s^{a-1}t^{b-1}.$
\end{remark}

\subsection{Trace map of positive words}

We can moreover compute the trace map for any \emph{positive} word in $F_2$,
namely for any word of the form $w=x^{a_1}y^{b_1}\dots x^{a_k}y^{b_k}$
where $$a_2,\dots,a_k,b_1,\dots,b_{k-1} \geq 1 \ \text{ and } a_1,b_k \geq 0.$$

We note that we can consider only words of the form $w=x^{a_1}y^{b_1}\dots x^{a_k}y^{b_k}$
where $a_1,b_1,\dots,a_k,b_k \geq 1$, and then we call $k$ the \emph{``length''} of this word.

Indeed, if $b_k=0$ then
$\tr(x^{a_1}y^{b_1}\dots x^{a_{k-1}}y^{b_{k-1}}x^{a_k}) = \tr(x^{a_1+a_k}y^{b_1}\dots x^{a_{k-1}}y^{b_{k-1}})$,
is the trace map of a positive word of length $k-1$.

\begin{theorem}\label{thm.gen.word}
Let $G=\SL(2,q)$ and let
$$w=x^{a_1}y^{b_1}\dots x^{a_k}y^{b_k}, \quad
a_1,b_1,\dots,a_k,b_k \geq 1.$$ Denote $s=\tr(x), u=\tr(xy),
t=\tr(y),$ and $A=\sum_{i=1}^k a_i$, $B=\sum_{i=1}^k b_i$.

Then $tr(w)=P(s,u,t)=\sum_{r=0}^{k} u^rp_r(s,t), $
where
\begin{itemize}
\item $p_k(s,t)=s^{A-k} t^{B-k} +\Phi(s,t)$
is a polynomial in $s,t$ and
$$\deg_s \Phi(s,t) \leq A-k, \ \deg_t \Phi(s,t) \leq B-k,$$
$$\deg_s \Phi(s,t)+\deg_t \Phi(s,t)<A+B-2k;$$
\item for all $r<k$,
$$\deg_s p_r(s,t)\leq A, \ \deg_t p_r(s,t)\leq B, \ \deg p_r(s,t)\leq A+B-2k. $$
\end{itemize}
\end{theorem}

\begin{proof}
The proof is by induction on $k$. The case $k=1$ was treated in
\lemref{lem.ab}.

We may always assume that $a_1\geq a_k.$ Let
$$w_1(x,y)=x^{a_{1}}y^{b_{1}} \dots x^{a_{k-1}}y^{b_{k-1}},$$
$$w_2 (x,y)=x^{a_k}y^{b_k},$$
$$w_3(x,y)=x^{a_{1}-a_k}y^{b_{1}} \dots x^{a_{k-1}}y^{b_{k-1}-b_k}.$$

Then, by \eqref{tr.formula},
\begin{equation}\label{eq.tr1}
\begin{aligned}
\tr(w) =& \tr(x^{a_{1}}y^{b_{1}} \dots x^{a_k}y^{b_k})\\
=&\tr(x^{a_{1}}y^{b_{1}} \dots x^{a_{k-1}}y^{b_{k-1}})\tr(x^{a_k}y^{b_k})-
\tr(x^{a_{1}-a_k}y^{b_{1}} \dots x^{a_{k-1}}y^{b_{k-1}-b_k})\\
=&\tr(w_1(x,y))\tr(w_2(x,y))-\tr(w_3(x,y)),\end{aligned}
\end{equation}

By the induction assumption we have
$$\tr(w_1)=P_1(s,u,t)=\sum_{r=0}^{k-1} u^r \tilde p_r(s,t),$$
$$\tr(w_2)=uf+h,$$
$$\tr(w_3)=Q(s,u,t),$$ and
\begin{itemize}
\item $\tilde p_{k-1}(s,t)=s^{A-a_k-k+1} t^{B-b_k-k+1} +\Phi_1(s,t); $
\item  $\deg_s \Phi_1(s,t) \leq A-a_k-k+1,$ \
$\deg_t \Phi_1(s,t) \leq B-b_k-k+1,$ and
$$\deg_s \Phi_1(s,t) + \deg_t \Phi_1(s,t) < A+B-a_k-b_k-2k+2;$$
\item $\deg_s f=a_k-1,\ \deg_t f=b_k-1, \  \deg_s h\leq a_k, \  \deg_t h\leq b_k \ \deg h \leq a_k+b_k-2;$
\item for all $r<k-1$, $ \deg_s \tilde p_{r}(s,t)\leq A-a_k, \ \deg_t \tilde p_{r}(s,t)\leq B-b_k,$
and $\deg\tilde p_{r}(s,t)\leq A+B-a_k-b_k-2k+2.$
\end{itemize}

We want to show that
\begin{equation}\label{eq.Q}
\deg_u  Q<k, \  \deg_s  Q\leq A, \  \deg_t  Q\leq  B,  \ \deg_s  Q+\deg_t  Q\leq A+B-2k.
\end{equation}
Since
$$\tr(w)=P(s,u,t)=P_1(s,u,t)(uf+h)-Q(s,u,t),$$
the theorem would follow from \eqref{eq.Q}.

\medskip

Consider the following cases.

\medskip

{\bf Case 1.}  $w_3(x,y)$ is a positive word.

Then its length is either $k-1,$ if both $a=a_{1}-a_k> 0$ and
$b=b_{k-1}-b_k> 0,$ or $k-2$ if $a=0$ or $b=0.$
Anyway, $\deg_u Q<k, $  \  $\deg_s  Q\leq A-2a_k\leq A,$  \  $\deg_t  Q\leq
B-2b_k\leq B,$ and $\deg_s  Q+\deg_t  Q\leq  A-2a_k+B-2b_k-2k+4\leq A+B-2k$,
as needed.

\medskip

{\bf Case 2.} $a_{1}-a_k>0,$   \ $b_{k-1}-b_k<0.$

Then
$$\tr(w_3)=Q(s,u,t)= \tr( w_4(x,y))\tr (y^{b_k-b_{k-1}})-\tr( w_5(x,y)), $$
where $$ w_4(x,y)=x^{a_1-a_k+a_{k-1}}y^{b_1}\dots x^{a_{k-2}}y^{b_{k-2}},$$
$$w_5(x,y)=x^{a_1-a_k}y^{b_1}\dots x^{a_{k-2}}y^{b_{k-2}}x^{a_{k-1}}y^{b_k-b_{k-1}}.$$
Thus, $w_4$ is a word of length $k-2$ and $w_5$ has length $k-1,$
and both are positive words. Let $ Q_1(s,u,t)=\tr( w_4(x,y)),$ \  $
Q_2(s,u,t)=\tr( w_5(x,y)).$

By the induction assumption, $\deg_u  Q_1=k-2, \ \deg_u  Q_2=k-1, $ and
$$\deg_s  Q_1\leq A-2a_k, \
\deg_t  Q_1\leq B-b_k-b_{k-1},$$
$$\deg_s  Q_1+\deg_t  Q_1\leq  A-2a_k+ B-b_k-b_{k-1}-2k+4;$$
$$\deg_s  Q_2\leq A-2a_k, \  \deg_t  Q_2\leq B-2b_{k-1},$$
$$\deg_s  Q_2+\deg_t  Q_2\leq A+B-2b_{k-1}-2a_k-2k+2\leq A+B-2k;$$

Moreover, $T(t)=\tr (y^{b_k-b_{k-1}})$ is a polynomial in $t$ of
degree $b_k-b_{k-1},$ thus
$$\deg_t  Q_1T(t)\leq B-b_k-b_{k-1}+ b_k-b_{k-1}= B-2b_{k-1},$$
$$ \deg_s  Q_1T +\deg_t Q_1T\leq A-2a_k+ B-2b_{k-1}-2k+4\leq A+B-2k.$$
Hence, for $Q= Q_1T(t)-Q_2$ condition \eqref{eq.Q} is valid.

\medskip

{\bf Case 3.} $a_{1}-a_k=0,$   \ $b_{k-1}-b_k<0.$

In this case  $w_3(x,y)=y^{b_1}\dots x^{a_{k-1}}y^{b_{k-1}-b_k}$ and
$$Q(s,u,t)=\tr (x^{a_2}\dots  x^{a_{k-1}}y^{b_{k-1}-b_k+b_1}).$$

If $b_{k-1}-b_k+b_1\geq 0$ then the word is positive of length $k-2$ and \eqref{eq.Q} is valid.

If $b_{k-1}-b_k+b_1<0$  then we perform the procedure described in {\bf Case 2}
for computing
$$Q(s,u,t)=\tr (x^{a_2}\dots  x^{a_{k-1}}y^{b_{k-1}-b_k+b_1})
= \tr( w_4(x,y))\tr (y^{b_k-b_{k-1}-b_1})-\tr( w_5(x,y)),$$ where
$$ w_4(x,y)=x^{a_2+a_{k-1}}y^{b_2}\dots x^{a_{k-2}}y^{b_{k-2}},$$
$$w_5(x,y)=x^{a_2}y^{b_2}\dots x^{a_{k-2}}y^{b_{k-2}}x^{a_{k-1}}y^{-b_{k-1}+b_k-b_1}.$$
The length of $w_4$ is $k-3$, and the length of $ w_5 $ is $k-2.$
Hence, $\deg_u Q=k-2.$ Moreover, $\deg_s Q\leq A - 2a_{k}\leq A$,
$\deg_t Q\leq B-2b_1-2b_{k-1} \leq B,$  and
$$\deg_s Q+\deg_t Q\leq A-2a_k+B-2b_1-2b_{k-1}-2k+6\leq A+B-2k .$$

Once more, we find out that conditions \eqref{eq.Q} are met by $Q.$
\end{proof}

\begin{remark}\label{non-positive-2}
Assume that $a_i,b_i\ne0$ but not necessarily positive. In view of
\remarkref{rem.non-positive-1}, for the word $w=x^{a_1}y^{b_1}\dots
x^{a_k}y^{b_k}$ equation \eqref{eq.tr1} implies that
\begin{equation}\label{eq.DEC}
\tr(w)=\sum\limits_{0}^k u^r G_r(s,t) \text{ and }
G_k(s,t)=\prod\limits_{i=1}^{k}f_{a_i,b_i}.
\end{equation}
\end{remark}

\section{Basic facts on the word $w(x,y)=x^ay^b$ and finite groups}\label{sect.grp}
In this section we present some elementary facts regarding the
surjectivity of the word map $w=x^ay^b$ on finite groups.

\begin{prop}\label{prop.a}
Let $G$ be a finite group and let $a$ be an integer. Then the word
map corresponding to $w=x^a$ is surjective on $G$ if and only if
$\gcd(|G|,a)=1$.
\end{prop}
\begin{proof}
Let $d=\gcd(|G|,a)$. If $d>1$ then there exists some prime $p$ which
divides both $a$ and $|G|$. Thus, $G$ contains some element $g \ne
id$ of order $p$, and so,
\[
    g^a = (g^p)^{a/p} = id.
\]
Hence the word map $w=x^a$ is not $1$ to $1$, and cannot be
surjective on $G$.

If $d=1$ then there exists some integer $l$ s.t. $l\cdot a \equiv 1
\pmod{|G|}$. Let $g \in G$ and take $x=g^l$, then
\[
    x^a = (g^l)^a = g^{l\cdot a} = g^{1} = g,
\]
as needed.
\end{proof}

\begin{prop}\label{prop.a.b.rel}
Let $G$ be a finite group and let $a,b$ be two relatively prime
integers. Then the word map corresponding to $w=x^ay^b$ is always
surjective on $G$.
\end{prop}
\begin{proof}
Since $a,b$ are relatively prime, there exist integers $k,l$ s.t.
$k\cdot a + l \cdot b = 1$. Let $g \in G$ and take $x=g^k$ and
$y=g^l$, then $x^ay^b = g^{k\cdot a}g^{l \cdot b} = g^{k\cdot a + l
\cdot b} = g^{1} = g$.
\end{proof}

\begin{prop}\label{prop.a.G.rel}
Let $G$ be a finite group and let $a,b$ be two integers.
If either $a$ or $b$ is relatively prime to $|G|$ then
the word map corresponding to $w=x^ay^b$ is surjective on $G$.
\end{prop}
\begin{proof}
Assume that $a$ is relatively prime to $|G|$. Then there exists some
integer $l$ s.t. $l\cdot a \equiv 1 \pmod{|G|}$. Then for every $g
\in G$, take $x=g^l$ and $y=id$. Thus,
\[
    x^ay^b = (g^l)^a \cdot id^b = g^{l\cdot a} \cdot id = g^{1} = g.
\]
\end{proof}

\begin{remark}\label{rem.a.b.mod}
Let $G$ be a finite group and let $w=x^ay^b$. We can always assume
that $0 \leq a,b < \exp(G)$. Moreover, if $G$ is of even order, we
can assume that $0 \leq a,b \leq \exp(G)/2$.

Indeed, let $a_1 = a \bmod {\exp(G)}$ and $b_1 = b \bmod \exp(G)$,
then for every $x,y \in G$, $x^ay^b = x^{a_1}y^{b_1}$. Hence,
$x^ay^b$ is surjective on $G$ if and only if $x^{a_1}y^{b_1}$ is
surjective on $G$.

If $\exp(G)$ is even, let $$a_2=\begin{cases} a_1 &\text{ if } a_1 \leq \exp(G)/2\\
\exp(G)-a_1 &\text{ if } a_1 > \exp(G)/2 \end{cases} \text{ and }
b_2=\begin{cases} b_1 &\text{ if } b_1 \leq \exp(G)/2\\
\exp(G)-b_1 &\text{ if } b_1 > \exp(G)/2 \end{cases}.$$
Then, for every $x,y \in G$, $x^{a_1}y^{b_1}=x^{\epsilon_1 a_2}y^{\epsilon_2 b_2}$,
where $\epsilon_1,\epsilon_2 \in \{\pm 1\}$, and
\begin{align*}
&\{z=x^{a_2}y^{b_2}:\ x,y \in G\} = \{z=x^{-a_2}y^{b_2}:\ x,y \in
G\}\\ =&\{z=x^{a_2}y^{-b_2}:\ x,y \in G\} = \{z=x^{-a_2}y^{-b_2}:\
x,y \in G\}.
\end{align*}
\end{remark}

\section{Properties of the groups $\SL(2,q)$ and $\PSL(2,q)$}\label{sect.prop.SL}
In this section we summarize some well-known properties of the
groups $\SL(2,q)$ and $\PSL(2,q)$ (see for example \cite{Do} and
\cite{Su}).

Let $q = p^e$, where $p$ is a prime number and $e \geq 1$. Recall
that $\GL(2,q)$ is the group of invertible $2 \times 2$ matrices
over the finite field with $q$ elements, which we denote by $\BF_q$,
and $\SL(2,q)$ is the subgroup of $\GL(2,q)$ comprising the matrices
with determinant $1$. Then $\PGL(2,q)$ and $\PSL(2,q)$ are the
quotients of $\GL(2,q)$ and $\SL(2,q)$ by their respective centers.
Also recall that $\PSL(2,q)$ is \emph{simple} for $q \neq 2,3$.

Denote $d=\gcd(2,q-1) = \begin{cases} 1 \ & \text{ if } q \text{ is  even} \\
2 \ & \text{ if } q \text{ is  odd}
\end{cases}.$

Then the orders of $\SL(2,q)$ and $\PSL(2,q)$ are $q(q-1)(q+1)$ and
$\frac{1}{d}q(q-1)(q+1)$ respectively, and their respective
exponents are $\frac{1}{d}p(q^2-1)$ and $\frac{1}{d^2}p(q^2-1)$.

One can classify the elements of $\SL(2,q)$ according to their
possible Jordan forms. The following Table \ref{table.sl2} lists the
three types of (non-central) elements, according to whether the
characteristic polynomial $P_t(\la):=\la^2 - t \la +1$ of the matrix
$A \in \SL(2,q)$ (where $t=\tr(A)$) has $0$, $1$ or $2$ distinct
roots in $\BF_q$.

\begin{table}[h]
\begin{tabular} {|c|c|c|c|c|c|}
\hline
element & roots       &  canonical form in       & order in & order in & conjugacy classes  \\
type    & of $P_t(\la)$ & $\SL(2,\overline{\BF}_p)$ & $\SL(2,q)$ & $\PSL(2,q)$ & in $\SL(2,q)$\\
\hline \hline

$id$ & $1$ root & $\begin{pmatrix} 1 & 0 \\ 0 & 1
\end{pmatrix}$ & $1$ & $1$ & one element \\
\hline

$-id$ & $1$ root & $\begin{pmatrix} -1 & 0 \\ 0 & -1
\end{pmatrix}$ & $d$ & $1$ & one element \\
\hline \hline

unipotent & $1$ root & $\begin{pmatrix} 1 & 1 \\ 0 & 1 \end{pmatrix}$ & $p$ & $p$ & $d$ conjugacy classes \\
 & & $t=2$ & & & each of size $\frac{q^2-1}{d}$ \\
\hline

 & $1$ root & $\begin{pmatrix} -1 & 1 \\ 0 & -1 \end{pmatrix}$ & $dp$ & $p$ & $d$ conjugacy classes \\
 & & $t=-2$ & & & each of size $\frac{q^2-1}{d}$ \\
\hline \hline

semisimple & $2$ roots & $\begin{pmatrix} \al & 0 \\ 0 & \al^{-1}\end{pmatrix}$ & divides & divides & for each $t$: \\
split      &        & where $\al \in \mathbb{F}_q^*$ & $q-1$ & $\frac{q-1}{d}$ & one conjugacy class \\
      &             & and $\al+\al^{-1}=t$ &  & & of size $q(q+1)$ \\
\hline \hline

semisimple & no roots & $\begin{pmatrix} \al & 0 \\ 0 & \al^q\end{pmatrix}$ & divides & divides & for each $t$: \\
non-split                 & & where $\al \in \mathbb{F}_{q^2}^* \setminus \mathbb{F}_q^*$ & $q+1$ & $\frac{q+1}{d}$ & one conjugacy class\\
         & & $\al^{q+1}=1$ &  & & of size $q(q-1)$ \\
 & & and $\al + \al^q = t$ & &  &\\
\hline \hline
\end{tabular}
\caption{Elements in the groups $\SL(2,q)$ and $\PSL(2,q)$.}
\label{table.sl2}
\end{table}

Table \ref{table.sl2}
shows that there is a deep connection between
\emph{traces} of elements in $\SL(2,q)$, their \emph{orders} and
their \emph{conjugacy classes}, as is expressed in the following
Lemmas.

\begin{Lemma}\label{lem.traces}${}$
\begin{itemize}
\item If $q$ is odd, then $x \in \SL(2,q)$ has order $4$ if and only if $\tr(x)=0$.
\item $x \in \SL(2,q)$ has order $3$ if and only if $\tr(x)=-1$.
\item If $p\geq 5$, then $x \in \SL(2,q)$ has order $6$ if and only if $\tr(x)=1$.
\end{itemize}

Moreover, for any $x \in \SL(2,q)$ satisfying $\tr(x)\ne 0,\pm 1, \pm 2$,
there exists some $y \in \SL(2,q)$ such that $\tr(x) \ne \tr(y)$, but the
orders of $x$ and $y$ are the same.
\end{Lemma}
\begin{proof}
If $m>2$ is an integer dividing $q-1$ then $\BF_q \setminus
\{0,1,-1\}$ contains $\phi(m)$ elements of order $m$, where $\phi$
denotes \emph{Euler's phi function}. Similarly, if $m>2$ divides
$q+1$ then $\BF_{q^2} \setminus \BF_q$ contains $\phi(m)$ elements
$\al$ of order $m$ satisfying $\al^{q+1}=1$.

Hence, if $m>2$ divides either $q-1$ or $q+1$, then
\[
    \# \{t \in \BF_q:\ t=\tr(x), \ x \in \SL(2,q),\ |x|=m \} =
    \frac{\phi(m)}{2}.
\]
The claim follows from the fact that $\phi(m) \geq 4$ if and only if
$m \neq 1,2,3,4,6$.
\end{proof}

\begin{Lemma}\label{lem.conj.unip}
Assume that $q$ is odd, take $\la \in \BF_{q^2} \setminus \BF_q$
satisfying $\la^2 \in \BF_q$, and let $g = \begin{pmatrix} \la & 0 \\
0 & \la^{-1} \end{pmatrix}$. Then for any $x \in \SL(2,q)$, $gxg^{-1} \in \SL(2,q)$, and
moreover,
\begin{itemize}
\item If $\tr(x)=2$ then exactly one of $x, gxg^{-1}$ is conjugate in $\SL(2,q)$ to
$\begin{pmatrix} 1 & 1 \\ 0 & 1 \end{pmatrix};$
\item If $\tr(x)=-2$ then exactly one of $x, gxg^{-1}$ is conjugate in $\SL(2,q)$ to
$\begin{pmatrix} -1 & 1 \\ 0 & -1 \end{pmatrix}.$
\end{itemize}
\end{Lemma}

\begin{proof}
Let $x = \begin{pmatrix} \al & \be \\ \ga & \de \end{pmatrix} \in
\SL(2,q)$, then $gxg^{-1} = \begin{pmatrix} \al & \be\la^2 \\
\ga\la^{-2} & \de \end{pmatrix} \in \SL(2,q)$.

Moreover, if $x = \begin{pmatrix} 1 & 1 \\ 0 & 1 \end{pmatrix}$ then
$gxg^{-1} = \begin{pmatrix} 1 & \la^{2} \\ 0 & 1 \end{pmatrix}$ is
not conjugate to $x$ in $\SL(2,q)$, since $\la^{2}$ is not a square
of some element in $\BF_q$.
\end{proof}

\begin{corollary}\label{cor.conj.tr}
Let $w \in F_2$ be some non-trivial word, let $z \ne \pm id$ be some
matrix in $\SL(2,q)$, and assume that $z$ can be written as
$z=w(x,y)$ for some $x,y \in \SL(2,q)$. Then for any matrix $\pm id
\ne z' \in \SL(2,q)$ with $\tr(z')=\tr(z)$ there exist $x',y' \in
\SL(2,q)$ such that $z'=w(x',y')$.
\end{corollary}
\begin{proof}
If $q$ is even, or if $q$ is odd and $\tr(z) \ne \pm 2$, then
necessarily $z'=hzh^{-1}$ for some $h \in \SL(2,q)$, so one can take
$x'=hxh^{-1}$ and $y'=hyh^{-1}$, and then
$$ w(x',y') =  w(hxh^{-1},hyh^{-1}) = hw(x,y)h^{-1} = hzh^{-1} =
z'.$$

Assume that $q$ is odd, and let $z' \in \SL(2,q)$ be some element
with $\tr(z')=2=\tr(z)$, then by \lemref{lem.conj.unip}, $z'$ is
either conjugate in $\SL(2,q)$ to $z$ or to $gzg^{-1}$ (where $g \in
\SL(2,q^2)$).

If $z'=hzh^{-1}$ for some $h \in \SL(2,q)$, take $x'=hxh^{-1}$ and
$y'=hyh^{-1}$, and then
$$ w(x',y') = w(hxh^{-1},hyh^{-1}) = hw(x,y)h^{-1} = hzh^{-1} =
z'.$$

If $z'=hgzg^{-1}h^{-1}$ for some $h \in \SL(2,q)$, take
$x'=hgxg^{-1}h^{-1}$ and $y'=hgyg^{-1}h^{-1}$, and then $x',y' \in
\SL(2,q)$ and moreover,
$$ w(x',y') = w(hgxg^{-1}h^{-1},hgyg^{-1}h^{-1}) = hgw(x,y)g^{-1}h^{-1} = hgzg^{-1}h^{-1} =
z'.$$

Similarly, if $\tr(z')=-2=\tr(z)$, then $z'=w(x',y')$ for some
$x',y' \in \SL(2,q)$.
\end{proof}

\section{Surjectivity of the trace map of $w(x,y)=x^ay^b$ on $\BF_q$}
\label{sect.surj.trace.map}
Recall that by \lemref{lem.ab}, the trace map of $w(x,y)=x^ay^b$ can
be written as
$$\tr  w(x,y)= u\cdot f_{a,b}(s,t)+h_{a,b}(s,t).$$

The following proposition shows that if neither $a$ nor $b$ is
divisible by the exponent of $\PSL(2,q)$, then the polynomial
$f_{a,b}(s,t)$ does not vanish identically on $\BA^2_{s,t}(\BF_q)$.

\begin{prop}\label{prop.summ.tr.surj}
Let $a,b\geq1$ and assume that neither $a$ nor $b$ is divisible by
$\frac{p(q^2-1)}{d^2}$. Then $f_{a,b}(s,t)$ does not vanish
identically on $\BA^2_{s,t}(\BF_q)$.

In particular, the following table summarizes the possible nine cases.
\begin{center}
\begin{tabular}{c||c|c|c|}
& $p \nmid a$ & $\frac{q-1}{d} \nmid a $ & $\frac{q+1}{d} \nmid a$ \\
\hline \hline
$p \nmid b$ & $f_{a,b}(2,2) \ne 0$ & $f_{a,b}(s_1,2) \ne 0$ & $f_{a,b}(s_2,2) \ne 0$ \\
\hline
$\frac{q-1}{d} \nmid b $ & $f_{a,b}(2,t_1) \ne 0$ & $f_{a,b}(s_1,t_1) \ne 0$ & $f_{a,b}(s_2,t_1) \ne 0$ \\
\hline
$\frac{q+1}{d} \nmid b $ & $f_{a,b}(2,t_2) \ne 0$ & $f_{a,b}(s_1,t_2) \ne 0$ & $f_{a,b}(s_2,t_2) \ne 0$\\
\hline
\end{tabular}
\end{center}
where:

$s_1 = \tr(x_1)$, $x_1$ is any element of order $q-1$;

$s_2 = \tr(x_2)$, $x_2$ is any element of order $q+1$;

$t_1 = \tr(y_1)$, $y_1$ is any element of order $q-1$;

$t_2 = \tr(y_2)$, $y_2$ is any element of order $q+1$.
\end{prop}

\begin{proof}
If $f_{a,b}(s,t)$ vanishes identically on $\BA^2_{s,t}(\BF_q)$,
then $\tr w(x,y)=h_{a,b}(s,t)$ does not depend on $u.$
We have to show that it is not the case for every $\BF_q.$
Take
$$x(\lambda, c)=\begin{pmatrix}\lambda & c \\ 0&\frac{1}{\lambda }\end{pmatrix},\quad
y(\mu, d)=\begin{pmatrix}\mu & 0 \\ d&\frac{1}{\mu}\end{pmatrix}.$$
Then for any $m,n$,
$$x(\lambda, u)^n=\begin{pmatrix}\lambda^n & c h_n(\lambda)\\0&\frac{1}{\lambda ^n}\end{pmatrix},\quad
y(\mu, v)^m=\begin{pmatrix}\mu^m & 0\\d h_n(\mu)&\frac{1}{\mu ^m}\end{pmatrix}.$$

\begin{Lemma}\label{lem.hn}
$$h_n(\zeta)=\frac{\zeta^{2n}-1}{\zeta^{n-1}(\zeta^{2}-1)}.$$
\end{Lemma}
\begin{proof}
We use induction on $n.$ For $n=1$ we have $h_n(\zeta)=1.$

Assume that for $n>1$ it is proved.
Then, by computing $x(\zeta,c)^{n+1}=x^nx$ and  $y(\zeta,d)^{n+1}=y^ny$ from the induction assumption
we obtain, respectively,
\begin{equation}\label{eq.hnn1}
h_{n+1}(\zeta)=\zeta^n+\frac{h_n(\zeta)}{\zeta}; \quad
h_{n+1}(\zeta)= \zeta h_n(\zeta)+\frac{1}{\zeta^n}.
\end{equation}
Both relation lead to the same result:
$$h_{n+1}(\zeta)=\zeta^n+\frac{\zeta^{2n}-1}{\zeta^{n}(\zeta^{2}-1)}=
\frac{\zeta^{2n}(\zeta^{2}-1)+\zeta^{2n}-1}{\zeta^{n}(\zeta^{2}-1)}=
\frac{\zeta^{2n+2}-1}{\zeta^{n}(\zeta^{2}-1)};$$

$$h_{n+1}(\zeta)=\zeta\frac{\zeta^{2n}-1}{\zeta^{n-1}(\zeta^{2}-1)}+\frac{1}{\zeta^n}=
\frac{\zeta^2(\zeta^{2n}-1)+(\zeta^{2}-1)}{\zeta^n(\zeta^{2}-1)}=
\frac{\zeta^{2n+2}-1}{\zeta^{n}(\zeta^{2}-1)}.$$
\end{proof}

Now, a direct computation shows that
$$\tr(x(\lambda, c)^a y(\mu, d)^b)=\lambda^a\mu^b+cd h_a(\lambda)h_b(\mu)+\frac{1}{\lambda^a\mu^b}.$$
We have to show that for every field $\BF_q$ there are $x(\lambda, c)$ and $y(\mu, d)$
such that
$$ h_a(\lambda)h_b(\mu)=\frac{(\lambda^{2a}-1)(\mu^{2b}-1)}
{(\lambda^{2}-1)(\mu^{2}-1)\lambda^{a-1}\mu^{b-1}}\ne 0.$$
Note that $h_n(1)=n.$

Let  $\al\in \BF_q$ be an element such that $\al^{q-1}=1, \al^{m}\ne 1 $ for any $m<q-1,$
and let $\be\in \BF_{q^2}\setminus\BF_q$ be an element satisfying that $\be^{q+1}=1, \be^{m}\ne 1$
for any $m<q+1.$

Since $\frac{p(q^2-1)}{d^2} \nmid a$ it follows that either $p$ or
$\frac{q-1}{d}$ or $\frac{q+1}{d}$ does not divide $a$. Similarly,
either $p$ or $\frac{q-1}{d}$ or $\frac{q+1}{d}$ does not divide $b$.

Thus, we need to consider nine cases, and in each case
we have to find  $x(\lambda, c)$ and $y(\mu, d)$  such that $h_a(\lambda)h_b(\mu)\ne 0.$
The following table shows that this is possible.
\begin{center}
\begin{tabular}{c||c|c|c|}
& $p \nmid a$ & $\frac{q-1}{d} \nmid a $ & $\frac{q+1}{d} \nmid a$ \\
\hline \hline
$p \nmid b$ & $\lambda=\mu=1$ & $\lambda=\al,\mu=1$& $\lambda=\be,\mu=1$ \\
\hline
$\frac{q-1}{d} \nmid b $ & $\lambda=1,\mu=\al$ & $\lambda=\al ,\mu=\al$& $\lambda=\be,\mu= \al$\\
\hline
$\frac{q+1}{d} \nmid b $ & $\lambda=1,\mu=\be$ & $\lambda=\al ,\mu=\be$& $\lambda=\be,\mu=\be$\\
\hline
\end{tabular}
\end{center}
\end{proof}

We can now deduce that if neither $a$ nor $b$ is divisible by the
exponent of $\PSL(2,q)$, then the trace map of the word
$w(x,y)=x^ay^b$ is surjective onto $\BF_q$.

\begin{corollary} \label{cor.tr.surj}
Let $a,b\geq1$ and assume that neither $a$ nor $b$ is divisible by
$\frac{p(q^2-1)}{d^2}$. Then every $\al \in \BF_q$ can be written as
$\al=\tr(x^ay^b)$ for some $x,y \in \SL(2,q)$.
\end{corollary}
\begin{proof}
According to \lemref{lem.ab} the trace of $w=x^ay^b$ can be written
as $$\tr(w)=u\cdot f(s,t)+h(s,t),$$ where $s=\tr(x), u=\tr(xy),
t=\tr(y).$ Namely, it is linear in $u$ and the coefficient of $u$ is
a non-trivial polynomial $f(s,t)$ in $s$ and $t$.

By \propref{prop.summ.tr.surj}, $f(s,t)$ does not vanish identically
on $\BA^2_{s,t}(\BF_q)$, and hence for every $\al \in \BF_q$ there
is a solution $(s,u,t) \in \BF_q^3$ to the equation $$u\cdot
f(s,t)+h(s,t) = \al.$$
\end{proof}

\section{Surjectivity of $w(x,y)=x^ay^b$ on $\SL(2,q)\setminus \{-id\}$ and $\PSL(2,q)$}
\label{sect.surj}

In this Section we prove Theorems~\ref{thm.surj.SL.even},
\ref{thm.surj.SL.odd} and \ref{thm.surj.PSL} on the surjectivity of
the word map $w(x,y)=x^ay^b$ on $\SL(2,q)\setminus \{-id\}$ and
$\PSL(2,q)$.

\medskip

By \remarkref{rem.a.b.mod}, throughout this section we can assume
that $1 \leq a,b \leq \frac{p(q^2-1)}{d^2}.$

The following two Corollaries follow from the general arguments
presented in \secref{sect.grp}.

\begin{corollary}\label{cor.a.b.d.surj}
If either $a$ is relatively prime to $\frac{p(q^2-1)}{d^2}$ or $b$
is relatively prime to $\frac{p(q^2-1)}{d^2}$, then the word
$w(x,y)=x^ay^b$ is surjective on $\SL(2,q)$, and hence on
$\PSL(2,q)$.
\end{corollary}
\begin{proof}
The claim follows immediately from \propref{prop.a.G.rel}.
\end{proof}

\begin{corollary}\label{cor.a.b.d.not.surj}
If either $a=\frac{p(q^2-1)}{d^2}$ and $b$ is not relatively prime
to $\frac{p(q^2-1)}{d^2}$; or $b=\frac{p(q^2-1)}{d^2}$ and $a$ is
not relatively prime to $\frac{p(q^2-1)}{d^2}$; then the word map
$w(x,y)=x^ay^b$ is \emph{not} surjective on $\PSL(2,q)$, and hence
not on $\SL(2,q) \setminus \{-id\}$.
\end{corollary}
\begin{proof}
The claim follows immediately from \propref{prop.a}.
\end{proof}

\begin{remark}
It follows that the only interesting cases to consider are when
$1\leq a,b<\frac{p(q^2-1)}{d^2}$, and both $a$ and $b$ are
\emph{not} relatively prime to $\frac{p(q^2-1)}{d^2}$.
\end{remark}

We can now deduce \propref{prop.abne2}, stating that if $w = x^ay^b$ is
non-degenerate with respect to $q$,
then any \emph{semisimple} element $z$ (namely, when $\tr(z)\ne \pm 2$) can be written as $z = x^ay^b$
for some $x,y \in \SL(2,q)$.

\begin{proof}[Proof of \propref{prop.abne2}]
Assume that $w = x^ay^b$ is non-degenerate with respect to $q$.
Without loss of generality, we may assume that $1 \leq a,b \leq \frac{p(q^2-1)}{d^2}.$
If $a,b<\frac{p(q^2-1)}{d^2}$ then the result immediately follows from
\corref{cor.tr.surj} and \secref{sect.prop.SL}. Otherwise, either
$a$ or $b$ is relatively prime to $\frac{p(q^2-1)}{d^2}$, and the
result follows from \corref{cor.a.b.d.surj}.
\end{proof}

Unfortunately, a similar result fails to hold when $z$ is
\emph{unipotent}, namely when $z \ne \pm id$ and $\tr(z)=\pm 2$.
This case will be discussed in detail in \secref{sect.non.surj},
where we shall prove the following two propositions.

\begin{prop} \label{prop.ab2}
Let $1 \leq a,b < \frac{p(q^2-1)}{d^2}$. Then the image of the word
map $w=x^ay^b$ contains any non-trivial element $z \in \SL(2,q)$
with $\tr(z)=2$, if and only if none of the following obstructions occurs:
\begin{enumerate}\renewcommand{\theenumi}{\it \roman{enumi}}
\item $q=2^e$, $e$ is odd and $a,b \in \{\frac{2(q^2-1)}{3},\frac{4(q^2-1)}{3}\}$;
\item $q \equiv 3 \bmod 4$ and $a=b=\frac{p(q^2-1)}{8}$;
\item $q \equiv 11 \bmod 12$ and $a=b=\frac{p(q^2-1)}{6}$;
\item $q \equiv 5 \bmod 12$ and $a,b \in \{\frac{p(q^2-1)}{6},
\frac{p(q^2-1)}{12}\}$.
\end{enumerate}
\end{prop}

\begin{prop} \label{prop.abmin2}
Assume that $q$ is odd and let $1 \leq a,b < \frac{p(q^2-1)}{4}$.
Then the image of the word map $w=x^ay^b$ contains any element $z
\ne -id$ with $\tr(z)=-2$, unless $q \equiv 3 \bmod 4$ and
$a=b=\frac{p(q^2-1)}{8}$.
\end{prop}


We can now prove the main theorems.

\begin{proof}[Proof of \thmref{thm.surj.SL.even}]
Let $q=2^e$. If $w=x^ay^b$ degenerates with respect to $q$, then by
\corref{cor.a.b.d.not.surj} the word map $w$ is not surjective on
$\SL(2,q)$. If obstruction {\it (i)} occurs then by
\propref{prop.ab2}{\it (i)}, the image of $w=x^ay^b$ does not
contain any non-trivial element $z \in \SL(2,q)$ with $\tr(z)=0$,
and hence it cannot be surjective on $\SL(2,q)$.

On the other hand, if $w=x^ay^b$ is non-degenerate with respect to
$q$ and obstruction {\it (i)} does not hold, then by
\propref{prop.abne2} and \propref{prop.ab2}{\it (i)}, any element $z
\in \SL(2,q)$ is in the image of the word map $w=x^ay^b$.
\end{proof}

\begin{proof}[Proof of \thmref{thm.surj.SL.odd}]
Let $q$ be an odd prime power. If $w=x^ay^b$ degenerates with
respect to $q$, then by \corref{cor.a.b.d.not.surj} the word map $w$
is not surjective on $\SL(2,q)\setminus \{-id\}$. If one of the
obstructions {\it (ii), (iii), (iv)} occurs then by
\propref{prop.ab2}, the image of $w=x^ay^b$ does not contain any
non-trivial element $z \in \SL(2,q)$ with $\tr(z)=2$, and hence it
cannot be surjective on $\SL(2,q) \setminus \{-id\}$.

On the other hand, if $w=x^ay^b$ is non-degenerate with respect to
$q$ and none of the obstructions {\it (ii), (iii), (iv)} holds, then
by \propref{prop.abne2}, \propref{prop.ab2} and
\propref{prop.abmin2}, any element $z \in \SL(2,q)\setminus \{-id\}$
is in the image of the word map $w=x^ay^b$.
\end{proof}

\begin{proof}[Proof of \thmref{thm.surj.PSL}]
Let $q$ be an odd prime power. Observe that the word map $w=x^ay^b$
is surjective on $\PSL(2,q)$ if and only if for every $z \in
\SL(2,q)$ either $z$ or $-z$ can be written as $x^ay^b$ for some
$x,y \in \SL(2,q)$.

If $w=x^ay^b$ degenerates with respect to $q$, then by
\corref{cor.a.b.d.not.surj} the word map $w$ is not surjective on
$\PSL(2,q)$. If obstruction {\it (ii)} occurs then by
\propref{prop.ab2} and \propref{prop.abmin2} the image of $w=x^ay^b$
does not contain any element $\pm id \ne z \in \SL(2,q)$ with
$\tr(z)=2$ or $\tr(z)=-2$, and hence it cannot be surjective on
$\PSL(2,q)$.

On the other hand, if $w=x^ay^b$ is non-degenerate with respect to
$q$ and obstruction {\it (ii)} does not hold, then by
\propref{prop.abne2}, \propref{prop.ab2} and \propref{prop.abmin2},
for any element $z \in \SL(2,q)$, either $z$ or $-z$ can be written
as $x^ay^b$ for some $x,y \in \SL(2,q)$, as needed.
\end{proof}

\begin{proof}[Proof of \corref{cor.Nab}]
If $q$ is odd, then by \thmref{thm.surj.PSL}, the word map
$w=x^ay^b$ is surjective on $\PSL(2,q)$ whenever
$$\frac{p(q^2-1)}{8} > \max \{a,b\}.$$

Thus, one has to prove that
\begin{equation}\label{eq.qodd1}
q > \sqrt{3a} \ \Longrightarrow \ \frac{p(q^2-1)}{8} > a
\end{equation} and
\begin{equation}\label{eq.qodd2}
|\PSL(2,q)| = \frac{q(q^2-1)}{2}> \frac{3\sqrt{3}}{2} a^{3/2} \
\Longrightarrow \ \frac{p(q^2-1)}{8} > a.
\end{equation}

Assume that $q > \sqrt{3a}$. Since $q\geq 3$ then
$q^2-1\geq\frac{8}{9}q^2.$ Therefore,
$$ \frac{p(q^2-1)}{8} \geq \frac{3(q^2-1)}{8}\geq \frac{3 \cdot 8 \cdot q^2}{9 \cdot 8} = \frac{q^2}{3} > a .$$
Moreover, the inequality
$$  \frac{3\sqrt{3}}{2} a^{3/2} < |\PSL(2,q)| = \frac{q(q^2-1)}{2} < \frac{q^3}{2} ,$$
implies that $q> \sqrt{3a}$.

This estimate is sharp. Indeed, if $a=p=q=3$, we have
$$q=3=\sqrt{3a}, \quad \frac{p(q^2-1)}{8}=3=a.$$

\medskip

If $q$ is even, then by \thmref{thm.surj.SL.even}, the word map $w=x^ay^b$
is surjective on $\PSL(2,q)$ whenever $$\frac{2(q^2-1)}{3} > \max \{a,b\}.$$

Thus, in this case one has to prove that
\begin{equation}\label{eq.qeven1}
 q > \sqrt{3a}\ \Longrightarrow \ \frac{2(q^2-1)}{3} >a
\end{equation} and
\begin{equation}\label{eq.qeven2}
 |\PSL(2,q)| = q(q^2-1)> \frac{3\sqrt{3}}{2} a^{3/2} \ \Longrightarrow \ \frac{2(q^2-1)}{3} >a.
\end{equation}

If $q>\sqrt{3a}$ then
$$\frac{2(q^2-1)}{3}\ge \frac{2(3a-1)}{3}=2a-\frac{2}{3}>a.$$

Let us prove \eqref{eq.qeven2}. If $q=2,$ then $q(q^2-1)=6\le
\frac{3\sqrt{3}}{2} a^{3/2}$ for any $a\geq 2.$ Hence, we may assume
that $q\geq 4,$ and then $\frac{2(q^2-1)}{3}\geq 10>3.$ It follows
that \eqref{eq.qeven2} is valid for $a=2,3$. On the other hand,
$$q(q^2-1)> \frac{3\sqrt{3}}{2} a^{3/2} \ \Longrightarrow \ q^3>\frac{3\sqrt{3}}{2} a^{3/2} $$
$$\Longrightarrow \ q^2>\frac{3a}{2^{\frac{2}{3}}}\ \Longrightarrow\
\frac{2(q^2-1)}{3}>2^{\frac{1}{3}}a-\frac{2}{3}\geq a$$ if $a \geq
3$.
\end{proof}

\section{Equidistribution of the word map $w(x,y)=x^ay^b$ on $\PSL(2,q)$}
\label{sect.equi}

The goal of this section is to prove \thmref{thm.equi}. We first
consider the case $\tilde G=\SL(2,q)$. In this case, the Theorem
follows from the following Proposition.

\begin{prop}\label{prop.equidist.sl2}
Denote $\tilde G=\SL(2,q)$, let $a,b \geq 1$ and let
$ w:\tilde G \times \tilde G \to \tilde G$ be the map
given by $w(x,y) = x^ay^b$. Then there are a subset $\tilde S
\subseteq \tilde G,$ and numbers $ A_1(a,b), A_2(a,b)$ such that:
\begin{enumerate}\renewcommand{\theenumi}{\it \roman{enumi}}
\item
$|\tilde S| = \bigl(1-\ep(q)\bigr)|\tilde G|,$  where
$0\leq\ep(q)\leq\frac{A_1(a,b)}{q};$
\item  $|M_g|= q^3(1+\de(q))$  for any $g\in \tilde  S,$
where $M_g=\{(x,y)\in \tilde G^2\ |  \ w(x,y)=g\}$ and $|\de(q)|\leq
\frac{A_2(a,b)}{q}. $\end{enumerate}
\end{prop}

\begin{proof}
We fix $a,b$ and maintain the notation of \secref{sect.trace.map}
omitting only the indices $a,b.$ Let $D=a+b-1.$

Consider the following commutative diagram of functions (between
finite sets):
\begin{equation}\label{diag.d1}
\xymatrix{ \tilde G \times \tilde G(\BF_q) \ar[r]^{\pi} \ar[d]^{w}
\ar[rd]^{\vp} & \BA^3_{s,u,t} (\BF_q) \ar[d]^{\psi}
\\ \tilde G  (\BF_q) \ar[r]^{\tau} & \BA^1_{s}(\BF_q)}
\end{equation}

In this diagram:
\begin{itemize}
\item $\BA^n_{x_1,x_2,\dots}$ denotes an $n-$dimensional affine space with coordinates $x_1,x_2,\dots;$
\item $\pi(x,y)=(\tr(x),\tr(xy),\tr(y))\in \BA^3_{s,u,t} ;$
\item $w(x,y)=x^ay^b \in \tilde G ;$
\item $\psi(s,u,t)=uf(s,t)+h(s,t) ;$
\item $\tau(x)=\tr(x) \in \BA^1_{s} ;$
\item $\vp(x,y)=\tr(w(x,y)) .$
\end{itemize}

By definition, $M_g=w^{-1}(g).$ Let $t \in \BF_q$. Denote:
\begin{itemize}
\item $N_t=\vp^{-1}(t)\subseteq \tilde G^2,$
\item $L_t=\psi^{-1}(t)\subseteq \BA^3_{s,u,t},$
\item $T_t=\tau^{-1}(t)\subseteq \tilde G,$
\item $p(s,u,t)=s^2+u^2+t^2-ust-4$,
\item $\nu_i(t), \ i=1,2,$ are solutions of the quadratic equation $x^2-tx+1=0,$
\item $\om_t^2=t^2-4.$
\end{itemize}

Note that $\nu_1(t) \ne \nu_2(t)$ if and only if $t \ne \pm 2$. For
odd $q$ the condition $\om_t \in \BF_q$ is equivalent to the
condition $\nu_{1,2} \in\BF_q$.

Recall that if $\pm 2\ne t\in \BF_q,$ then all the elements in $T_t$
are conjugate (see \secref{sect.prop.SL}). Thus, if $\tr(g)=t,$ then
$|M_g|=\frac{|N_t|}{|T_t|}$.
From Table~\ref{table.sl2} in \secref{sect.prop.SL} we deduce that
\begin{equation}\label{eq.Do}
|T_t|=q^2(1+\de_1(t)),
\end{equation}
where
$$
\de_1(t)=\begin{cases}0 & \text{ if } \om_t=0\\
\fc{1}{q} & \text{ if } \om_t\ne 0 \ and \ \nu_{1,2}(t) \in\BF_q \\
\fc{-1}{q}  & \text{ if } \om_t\ne 0 \ and \ \nu_{1,2}(t)  \not\in\BF_q
\end{cases}.$$
Hence, in all the cases above, $|\de_1|\leq\frac{1}{q}$.

\medskip

We divide the proof into three steps.

\subsection*{\bf Step 1} \emph{Fibers of $\pi.$}

\begin{prop}\label{prop.pi}
\begin{enumerate}\renewcommand{\theenumi}{\alph{enumi}}
\item If $t^2\ne 4,$ then

$|\pi^{-1}(s,u,t)(\BF_q)| \
\begin{cases}
= q^3(1+\de_2(t)) &\text{ if } p(s,u,t)\ne 0 \\
\leq 2q^3(1 +\frac{1}{q}) &\text{ if } p(s,u,t)= 0
\end{cases},$

where $|\de_2|\leq \frac{3}{q}$ for every $(s,u,t).$

\item
$|\pi^{-1}(s,u,\pm 2)(\BF_q)| \
\begin{cases}
= q^3-q &\text{ if } p(s,u,2)\ne 0 \\
\leq 2q^3(1+\de_1(s)) &\text{ if } p(s,u,2)= 0
\end{cases}.$
\end{enumerate}
\end{prop}

The proof of this Proposition follows from the next two Lemmas.

For a fixed $y \in \tilde G$ with $\tr(y)=t$, let
$$K_{s,u}(y)=\{x\in \tilde G \  | \ \pi(x,y)=(s,u,t)\}.$$

\begin{Lemma}\label{lem.yt}
Let $t^2\ne 4$ and $y_t=\begin{pmatrix} t& 1 \\ -1& 0
\end{pmatrix}$.
Then
\begin{equation}\label{eq.K1} |K_{s,u}(y_t)(\BF_q)|=\begin{cases}
q\pm 1 \ & \text{ if } p(s,u,t)\ne 0 \\
1 \ & \text{ if } p(s,u,t)= 0,   \ \nu_{1,2}(t) \not \in \BF_q\\
2q-1 \ & \text{ if } p(s,u,t)= 0, \    \nu_{1,2}(t) \in \BF_q
\end{cases}.
\end{equation}
\end{Lemma}

\begin{proof}
If $\om_t\ne 0$ then $K_{s,u}(y_t)$ consists of the matrices
\begin{equation*}\label{eq.f4}
x=\begin{pmatrix} \al & \be \\ u+\be-\al t & s-\al
\end{pmatrix},\end{equation*}
such that
\begin{equation}\label{eq.p201} \al(s-\al)-(u+\be-\al t)\be=1.\end{equation}

Assume that $q$ is odd. Then \eqref{eq.p201} is equivalent to
$$\left( \fc{st-2u}{2}+\fc{\om_t^2
\be}{2}\right)^2- \om_t^2\left(\al-\fc{\be
t}{2}-\fc{s}{2}\right)^2-p(s,u,t)=0. $$

Thus, if $p(s,u,t)\ne 0$ then $K_{s,u}(y_t)$ is a non-degenerate
conic; whereas if $p(s,u,t)=0$ then $K_{s,u}(y_t)$ is a pair of
intersecting straight lines, which are not defined over $\BF_q$ if
$\om_t \not \in \BF_q$, or defined over $\BF_q$ if $\om_t \in
\BF_q$.

Assume that $q=2^e.$ Substituting $\tilde \al=\al+\frac{u}{t},$
\ $\tilde \be=\be+\frac{s}{t},$ we reduce \eqref{eq.p201} to
\begin{equation}\label{eq.p202}
\tilde \al^2+\tilde \be^2+t\tilde \al\tilde \be=\frac{p(s,u,t)}{t^2}.
\end{equation}
Thus, since  $t\ne 0,$ we have a conic for  $p(s,u,t)\ne 0.$
Hence, $|K_{s,u}(y_t)(\BF_q)|=q\pm 1 $ if $p(s,u,t)\ne 0.$

If $p(s,u,t)= 0$ then
\eqref{eq.p202} provides
\begin{equation}\label{eq.p203}
(\nu_1(t)\tilde \al+\tilde \be)(\nu_2(t)\tilde \al+\tilde \be)=0\end{equation}
Thus, if $\nu_{1,2}(t)\in\BF_q$ we have two intersecting lines, whereas
if $\nu_{1,2}(t)\not \in\BF_q$ we have precisely one point $(0,0).$
\end{proof}

\begin{Lemma}\label{lem.vla} Let $t=2,$  $\la\ne 0 ,$ and
$v_{\la}=\begin{pmatrix} 1& \la \\ 0& 1 \end{pmatrix}$. Then
\begin{equation}\label{eq.K2} |K_{s,u}(v_\la)(\BF_q)|=\begin{cases}
q \ & \text{ if } p(s,u,2)\ne 0 \\
0 \ & \text{ if } p(s,u,2)= 0,  \ \nu_{1,2}(s) \not \in \BF_q\\
2q \ or \  q\ & \text{ if } p(s,u,2)= 0,  \   \nu_{1,2}(s) \in \BF_q
\end{cases}.\end{equation}
\end{Lemma}

\begin{proof}
If $\la\ne 0$ then $K_{s,u}(v_\la)$ consists of the matrices
\begin{equation*}\label{f8}
x=\begin{pmatrix} \al & \be \\ \frac{u-s}{\la} & s-\al
\end{pmatrix},\end{equation*}
such that $$\al(s-\al)-\be\frac{u-s}{\la}=1.$$

Thus if $p(s,u,2)=(s-u)^2\ne 0$ then we have $q$ points; whereas if
$p(s,u,2)=(s-u)^2= 0,$ then we have
\begin{itemize}
\item either two disjoint lines, if $\nu_{1,2}(s)\in \BF_q,\ \nu_1\ne\nu_2$
(for odd $q$ it means that $w_s\ne 0,w_s\in \BF_q)$;
\item or one line, if $\nu_{1,2}(s)\in \BF_q, \nu_1=\nu_2$ (for odd $q$ it means that $w_s=0$);
\item no points, if  $\nu_{1,2}(s)\not \in \BF_q,$  (for odd $q$ it means that $w_s\not\in \BF_q)$.
\end{itemize}
\end{proof}

\begin{proof}[Proof of \propref{prop.pi}] ${ }$

\begin{enumerate}\renewcommand{\theenumi}{\alph{enumi}}
\item Indeed, $|\pi^{-1}(s,u,t)(\BF_q)|=|K_{s,u}(y_t)(\BF_q)|\cdot|T_t|, $
thus the claim follows from \lemref{lem.yt} and
Equation~\eqref{eq.Do}. Moreover,
$$|\de_2(t)|\leq |\de_1(t)|+\frac{1}{q}+|\de_1(t)|\frac{1}{q}\leq \frac{3}{q}.$$

\item There are two or three conjugacy classes of matrices $y$ with $\tr(y)=2.$ If
$u\ne s,$ then $y\ne id,$ thus there are $q^2-1$ different matrices
$y$ to consider, and according to Equation~\eqref{eq.K2},
$$|\pi^{-1}(s,u,\pm 2)(\BF_q)|=(q^2-1)q.$$

If $p(s,u,2)= 0,$ i.e. $s=u,$ then summation over the classes
yields
$$|\pi^{-1}(s,u,\pm 2)(\BF_q)|\leq 2q(q^2-1) + q^2(1+\de_1(s))\leq 2q^3\bigl(1+\frac{1}{q}\bigr).$$

For $t=-2$ the proof is similar.
\end{enumerate}
\end{proof}

\subsection*{\bf Step 2} \emph{Definition of the set $\tilde S.$}

Let
$$A=\{(s,t)\in \BA^2_{s,t} \ |  \ f(s,t)=0\},$$
$$B_\zeta=\{(s,t)\in \BA^2_{s,t} \ | \  h(s,t)=\zeta \},$$
$$C=\{(s,u,t)\in \BA^3_{s,u,t} \ | \  p(s,u,t)=0\}.$$

Note that $C$ is absolutely irreducible for every field $\BF_q$.
We first define the set $\Sigma\subset \BF_q$ by the
following rules.

\begin{itemize}
\item \emph{Rule 1.} Assume that there exists $\zeta \in \BF_q$ satisfying
$p\bigl(s, \frac{\zeta-h(s,t)}{f(s,t)},t\bigr)\equiv 0$ on $C.$ Then
$\zeta \in \Sigma.$

\item \emph{Rule 2.} Assume that there is an irreducible (over
$\overline{\BF}_q$) component $A'\subseteq A$ and $\zeta\in \BF_q$
such that $h(s,t)\equiv \zeta$ on $A'.$  Then $\zeta \in \Sigma.$
Note that since $\deg A\leq D$ (by \lemref{lem.ab}) there are at most
$D$ such numbers.

\item \emph{Rule 3.} $2\in\Sigma$ and $-2\in\Sigma.$
\end{itemize}

\begin{remark}
By the above construction, $\Sigma$ contains all the values $\zeta$
such that $A\cap B_\zeta$ contains a curve or $L_\zeta\cap C$ is not
a curve.
\end{remark}

Now, we can define the sets $\tilde T=\tau^{-1}(\Sigma)$ and $\tilde
S=\tilde G\setminus \tilde T.$

\begin{Lemma}\label{lem.def.S}
$|\tilde S|=|\tilde G|(1-\varepsilon(q)),$ where
$|\varepsilon(q)|\leq\frac{3+D}{q-1}$.
\end{Lemma}
\begin{proof}
Indeed, by construction, $|\Sigma|\leq(3+D).$ Thus by \eqref{eq.Do},
$|\tilde T|\leq(3+D)q^2\bigl(1+\fc{1}{q}\bigr).$ Hence,
$$|\tilde S|=|\tilde G|-|\tilde T|=|\tilde G|(1-\varepsilon(q)),$$
where $$|\varepsilon(q)|\leq
\frac{(3+D)q^2(1+\fc{1}{q})}{q^3-q}=\frac{3+D}{q-1}.$$
\end{proof}

\subsection*{\bf Step 3} \emph{Estimation of $|M_g|.$}

Let $\zeta \in \BF_q \setminus \Sigma.$ Then
$L_\zeta=\psi^{-1}(\zeta)=Y_\zeta \cup R_\zeta \cup Q_\zeta,$ where
\begin{align*}
Y_\zeta&=\{(s,u,t)\in L_\zeta \ | \  p(s,u,t)\ne 0, f(s,t)\ne 0 \},
\\ R_\zeta&=\{(s,u,t)\in L_\zeta \ | \  p(s,u,t)=0 , f(s,t)\ne 0
\},\\ Q_\zeta&=\{(s,u,t)\in L_\zeta \ |  \  f(s,t)= 0 \}.
\end{align*}

In the estimation of the sizes of the above sets, we will use the
following fact, which is the case $n=1$ of \cite[Proposition
12.1]{GL}.

\begin{claim}
Let $X\subseteq \BP^N$ be a projective curve in the projective space
$\BP^N$ of degree $D$ defined over $\BF_q.$ Then
\begin{equation}\label{eq.lachaud}
|X(\BF_q)|\leq D(q+1).
\end{equation}
\end{claim}

\begin{Lemma}\label{lem.Qz}
$|Q_\zeta(\BF_q)|\leq D^2q.$
\end{Lemma}
\begin{proof}
If $f(s,t)= 0$ and $(s,u,t)\in L_\zeta,$ then $h(s,t)=\zeta,$
i.e.$(s,t)\in A \cap B_\zeta .$ Since $\zeta \not\in\Sigma,$  the
set $A\cap B_\zeta$ is finite. Both curves have degree at most $D$
(by \lemref{lem.ab}), hence, by B\'{e}zout's Theorem, $|A\cap
B_\zeta|\leq D^2.$ On the other hand there is no restriction on the
value of $u.$ Hence, $|Q_\zeta(\BF_q)|\leq D^2q.$
\end{proof}

\begin{Lemma}\label{lem.Rz}
$|R_\zeta(\BF_q)|\leq 3D(q+1).$
\end{Lemma}
\begin{proof}
Indeed, $R_\zeta(\BF_q)=\{p(s,u,t)=0, uf(s,t)+h(s,t)=\zeta\}$ is a
curve, since $\zeta\not\in\Sigma.$  By B\'{e}zout's Theorem,
$\deg(R_\zeta)\leq 3D.$ Hence, according to Equation
\eqref{eq.lachaud}, $|R_\zeta(\BF_q)|\leq 3D(q+1).$
\end{proof}

\begin{Lemma}\label{lem.Yz}
$|Y_\zeta(\BF_q)|=q^2(1+\de_3(q))$ where $|\de_3(\zeta)|\leq\frac
{2D}{q}.$
\end{Lemma}
\begin{proof}
Since $\deg A\leq D$, by Equation \eqref{eq.lachaud}, $|A(\BF_q)|\leq
D(q+1).$ Hence
$$q^2-|(\BA^2_{s,t}\setminus A)(\BF_q)|\leq D(q+1).$$
For every point $(s,t)\in (\BA^2_{s,t}\setminus A)(\BF_q)$ there is
precisely one point $(s, \frac{\zeta-h(s,t)}{f(s,t)},t)\in Y_\zeta.$
Thus,
$$|Y_\zeta(\BF_q)|=|(\BA^2_{s,t}\setminus A)(\BF_q)|=q^2(1+\de_3(q)),$$
and $$|\de_3(\zeta)|\leq \frac {D(q+1)}{q^2}\leq\frac{2D}{q}.$$
\end{proof}

We can now estimate $|M_g|=\frac{|N_\zeta|}{|T_\zeta|}.$ By
\propref{prop.pi}, we have
$$|\pi^{-1}(Y_\zeta)(\BF_q)|=q^5(1+\de_4),$$
where $$|\de_4|\leq \frac {2D}{q}+\frac {3}{q}+\frac {2D}{q}\frac
{3}{q}\leq \frac{8D+3}{q},$$ and
$$|\pi^{-1}(R_\zeta\cup Q_\zeta)(\BF_q)|\leq \bigl(D^2q+3D(q+1)\bigr)2q^3\bigl(1+\frac {1}{q}\bigr)\leq
2(D^2+6D)q^4\bigl(1+\frac {1}{q}\bigr).$$

Therefore
\begin{align*}
||N_\zeta(\BF_q)|-q^5(1+\de_4)|
&=|\pi^{-1}(L_\zeta)(\BF_q)|-|\pi^{-1}(Y_\zeta)(\BF_q)| \\
&\leq 2(D^2+6D)q^4\bigl(1+\frac {1}{q}\bigr),
\end{align*}
implying that
$$|N_\zeta(\BF_q)|=q^5(1+\de_5),$$
where $$\de_5\leq\frac{8D+3}{q}
+\frac{2(D^2+6D)(1+\frac{1}{q})}{q}\leq\frac{4D^2+32D+3}{q}.$$

We conclude that $$|M_g|=\frac
{|N_\zeta|}{|T_\zeta|}=\frac{q^5(1+\de_5)}{q^2(1+\de_1)}=q^3(1+\de_6(\zeta)),$$
where $$|\de_6(\zeta)|\leq\frac{4D^2+32D+3}{q}+\frac{2}{q} = \frac{4D^2+32D+5}{q}.$$

\medskip

For completing the proof of \propref{prop.equidist.sl2} it is
sufficient to take
$$A_1(a,b)=2(3+D)\ge\frac{(3+D)q}{q-1},$$
and $$A_2(a,b)=4D^2+32D+5.$$
\end{proof}

We can now prove \thmref{thm.equi} for the group $G=\PSL(2,q)$.

\begin{proof}[Proof of \thmref{thm.equi} for $G=\PSL(2,q)$]
Assume that $q$ is odd, denote $G=\PSL(2,q)$ and consider the commutative diagram
\begin{equation}\label{diag.d2}
\xymatrix{ \tilde G \times \tilde G(\BF_q) \ar[r]^{w_1}
\ar[d]^{\rho'} \ar[rd]^{\varkappa} & \tilde G(\BF_q) \ar[d]^{\rho}\\
G\times G(\BF_q) \ar[r]^{w_2} & G(\BF_q)}
\end{equation}
where
\begin{itemize}
\item
$\rho:\tilde G\to G$ is the natural projection $ \tilde G\to \tilde
G/Z(\tilde{G});$
\item
$\rho':\tilde G\times \tilde G\to G\times G$ is the projection
induced by $\rho;$
\item $w_1,w_2$ correspond to the map $(x,y)\to x^ay^b$ on $\tilde G\times \tilde G$ and on $G\times G$
respectively.
\end{itemize}

Define $S=\rho(\tilde S).$ Since for any $z\in G$, $\rho^{-1}(z)$
contains precisely two elements of $\tilde G,$ then
\propref{prop.equidist.sl2} implies that
$$|S|=\frac{|\tilde G|(1-\varepsilon(q))}{2}=|G|(1-\varepsilon(q)).$$

Take $z\in S,$ then $\rho^{-1}(z)=\{z_1,z_2\},$ and denote $H_z=w_2^{-1}(z).$
Let $y\in \tilde G$ and denote $M_y=w_1^{-1}(y)$. Then $M_{z_1}\cup
M_{z_2}=\rho^{-1}(w_2^{-1}(z))=\rho^{-1}(H_z).$

Thus,
$$|H_z|=\frac{|M_{z_1}|+|M_{z_2}|}{4}=\frac{2q^3(1+\de(q))}{4}=|G|(1+\de'(q)),$$
where $|\de'(q)|\to 0$ when $q\to \infty.$
\end{proof}

\section{Non-surjectivity of some words $w(x,y)=x^ay^b$}\label{sect.non.surj}

In this section we prove Propositions~\ref{prop.ab2}
and~\ref{prop.abmin2}, and in particular, we show that there are
certain fields $\BF_q$ and positive integers $a,b$ such that the
trace map corresponding to the word $w=x^ay^b$ is surjective on
$\BF_q$, by \corref{cor.tr.surj}, however, the word $w=x^ay^b$
itself is not surjective on $\SL(2,q)$ (or $\PSL(2,q)$), since the
image of $w$ does not contain $-id$ or unipotent elements, yielding
the obstructions described in Definition~\ref{def.obs}.

\subsection{Proof of Propositions~\ref{prop.ab2} and~\ref{prop.abmin2}}

\begin{proof}[Proof of~\propref{prop.ab2}]
Since $\frac{p(q^2-1)}{d^2} \nmid a$ it follows that either $p$ or
$\frac{q-1}{d}$ or $\frac{q+1}{d}$ does not divide $a$. Similarly,
either $p$ or $\frac{q-1}{d}$ or $\frac{q+1}{d}$ does not divide
$b$.

We may consider the following four cases:

\medskip

{\bf Case 1:} $p\nmid a$.

Take $x= \begin{pmatrix} 1 & \frac{1}{a}
\\0 & 1\end{pmatrix}$ and $y=id$. Then $x^ay^b=x^a=\begin{pmatrix} 1
& 1 \\0 & 1\end{pmatrix}$, and the claim follows from \corref{cor.conj.tr}.

\medskip

For the other cases, it is sufficient to find some $x,y \in
\SL(2,q)$ with $s=\tr(x), t=\tr(y)$ satisfying:
\begin{equation}\label{eq.tr.2}
\quad f_{a,b}(s,t) \ne 0 \text{ and } \tr(x^a) \ne \tr(y^b).
\end{equation}

Indeed, if $f_{a,b}(s,t) \ne 0$ then one can find some $u \in \BF_q$
such that
$$\tr(w)=u\cdot f(s,t)+h(s,t) = 2.$$

Hence, there exist some matrices $x_1,y_1 \in \SL(2,q)$ satisfying
$\tr(x_1)=s, \tr(y_1)=t$, $\tr(x_1y_1)=u$ and $\tr(x_1^ay_1^b)=2$.
Moreover, $x_1^ay_1^b \ne id$ since
$$\tr(x_1^a) = \tr(x^a) \ne \tr(y^b) = \tr(y_1^b).$$
Therefore, by \corref{cor.conj.tr}, for any $z \ne id$ with
$\tr(z)=2$ there exist some $x_2,y_2 \in \SL(2,q)$ such that
$z=x_2^a y_2^b$ as needed.

\medskip

{\bf Case 2:} $\frac{q-1}{d}\nmid a$ and $\frac{q+1}{d} \nmid b$.

Let $x$ and $y$ be two matrices of orders $q-1$ and $q+1$
respectively, and let $s=\tr(x)$ and $t=\tr(y)$.
According to the table in \propref{prop.summ.tr.surj},
$f_{a,b}(s,t) \ne 0$. Moreover, since $x$
is a split element while $y$ is a non-split element, and since $x^a
\ne \pm id$, $y^b \ne \pm id$, then $\tr(x^a) \ne \tr(x^b)$,
implying \eqref{eq.tr.2}.

\medskip
{\bf Case 3:} $\frac{q-1}{d}\nmid a$ and $\frac{q-1}{d} \nmid b$.

Let $x,y$ be some matrices of order $q-1$ and note that $x^a \ne \pm
id$ and $y^b \ne \pm id$. Observe that unless either $|x^a|=|y^b|=3$
or $|x^a|=|y^b|=4$, for all elements $x,y$ of order $q-1$, one can
find two matrices $x,y$ of order $q-1$ satisfying $\tr(x^a) \ne
\tr(y^b)$ (see \lemref{lem.traces}). Let $s=\tr(x)$ and $t=\tr(y)$.
According to the table in \propref{prop.summ.tr.surj},
$f_{a,b}(s,t) \ne 0$, and so \eqref{eq.tr.2} holds.

Since $a,b \leq \frac{p(q^2-1)}{2}$ the only cases left to consider
are cases $(i),(iii),(v)$ of \remarkref{rem.obs.3.4}. In
Proposition~\ref{prop.good.xaya} we will show that in all these
cases the image of $w=x^ay^b$ contains every non-trivial element $z
\in \SL(2,q)$ with $\tr(z)=2$.

\medskip
{\bf Case 4:} $\frac{q+1}{d}\nmid a$ and $\frac{q+1}{d} \nmid b$.

Let $x,y$ be some matrices of order $q+1$ and note that $x^a \ne \pm
id$ and $y^b \ne \pm id$. Similarly to Case 3, observe that unless
either $|x^a|=|y^b|=3$ or $|x^a|=|y^b|=4$, for all elements $x,y$ of
order $q+1$, one can find two matrices $x,y$ of order $q+1$
satisfying $\tr(x^a) \ne \tr(y^b)$ (see \lemref{lem.traces}). Let
$s=\tr(x)$ and $t=\tr(y)$. According to the table in \propref{prop.summ.tr.surj},
$f_{a,b}(s,t) \ne 0$, and so \eqref{eq.tr.2} holds.

Since $a,b \leq \frac{p(q^2-1)}{2}$ the only cases left to consider
are cases $(ii),(iv),(vi)$ of \remarkref{rem.obs.3.4}. In
Proposition~\ref{prop.obs.xaya} we will show that in all these cases
the image of $w=x^ay^b$ contains no non-trivial element $z \in
\SL(2,q)$ with $\tr(z)=2$, yielding the obstructions given in the
proposition.
\end{proof}

\begin{proof}[Proof of~\propref{prop.abmin2}]
Since $\frac{p(q^2-1)}{4} \nmid a$ it follows that either $p$ or
$\frac{q-1}{2}$ or $\frac{q+1}{2}$ does not divide $a$. Similarly,
either $p$ or $\frac{q-1}{2}$ or $\frac{q+1}{2}$ does not divide
$b$.

We may consider the following six cases:

\medskip

{\bf Case 1:} $p\nmid a$ and $p \nmid b$.

Take: $x= \begin{pmatrix} 1 & -2/a \\ 0 & 1 \end{pmatrix}, \quad
y= \begin{pmatrix} 1 & 0 \\ 2/b & 1 \end{pmatrix}.$

Then
$x^a= \begin{pmatrix} 1 & -2 \\ 0 & 1 \end{pmatrix}, \quad
y^b= \begin{pmatrix} 1 & 0 \\ 2 & 1 \end{pmatrix},$
and so
$z= x^ay^b = \begin{pmatrix} -3 & -2 \\ 2 & 1 \end{pmatrix} \ne -id,$
satisfies that $\tr(z)=-2$, and the result follows from \corref{cor.conj.tr}.

\medskip

For the other cases, it is sufficient to find some $x,y \in
\SL(2,q)$ with $s=\tr(x), t=\tr(y)$ satisfying:
\begin{equation}\label{eq.tr.min.2}
\quad f_{a,b}(s,t) \ne 0 \text{ and } \tr(x^a) \ne -\tr(y^b).
\end{equation}

Indeed, if $f_{a,b}(s,t) \ne 0$ then one can find some $u \in \BF_q$
such that
$$\tr(w)=u\cdot f(s,t)+h(s,t) = -2.$$

Hence, there exist some matrices $x_1,y_1 \in \SL(2,q)$ satisfying
$\tr(x_1)=s, \tr(y_1)=t$, $\tr(x_1y_1)=u$ and $\tr(x_1^ay_1^b)=-2$.
Moreover, $x_1^ay_1^b \ne -id$ since
$$\tr(x_1^a) = \tr(x^a) \ne -\tr(y^b) = -\tr(y_1^b).$$
Therefore, by \corref{cor.conj.tr}, for any $z \ne -id$ with $\tr(z)
\ne -2$ there exist some $x_2,y_2 \in \SL(2,q)$ such that $z=x_2^a
y_2^b$ as needed.

\medskip
{\bf Case 2:} $p\nmid a$ and $\frac{q-1}{2} \nmid b$.

Let $y$ be some matrix of order $q-1$ and let $t=\tr(y)$. According
to the table in \propref{prop.summ.tr.surj}, $f_{a,b}(2,t) \ne 0$.
Moreover, since $y$ is a split element and $y^b \ne \pm id$ then
$\tr(y^b) \ne \pm 2$, implying \eqref{eq.tr.min.2}.

\medskip
{\bf Case 3:} $p\nmid a$ and $\frac{q+1}{2} \nmid b$.

The proof is the same as in Case 2.

\medskip
{\bf Case 4:} $\frac{q-1}{2} \nmid a$ and $\frac{q+1}{2} \nmid b$.

Let $x$ and $y$ be two matrices of orders $q-1$ and $q+1$
respectively, and let $s=\tr(x)$ and $t=\tr(y)$.
According to the table in \propref{prop.summ.tr.surj},
$f_{a,b}(s,t) \ne 0$. Moreover, since $x$
is a split element while $y$ is a non-split element, and since $x^a
\ne \pm id$, $y^b \ne \pm id$, then $\tr(x^a) \ne \pm \tr(y^b)$,
implying \eqref{eq.tr.min.2}.

\medskip
{\bf Case 5:} $\frac{q-1}{2} \nmid a$ and $\frac{q-1}{2} \nmid b$.

Let $x,y$ be some elements of order $q-1$ and note that $x^a \ne \pm
id$ and $y^b \ne \pm id$. Observe that unless $\tr(x^a)=\tr(y^b)=0$
for all elements $x,y$ of order $q-1$, one can find two matrices
$x,y$ of order $q-1$ satisfying $\tr(x^a) \ne -\tr(y^b)$ (see
\lemref{lem.traces}). Let $s=\tr(x)$ and $t=\tr(y)$.
According to the table in \propref{prop.summ.tr.surj},
$f_{a,b}(s,t) \ne 0$, and so
\eqref{eq.tr.min.2} holds.

Now, the only case left is $q \equiv 1 \bmod 4$ and
$a=b=\frac{p(q^2-1)}{8}$ (see \remarkref{rem.obs.3.4}$(iii)$). In
\propref{prop.good.xaya} we will show that in this case the image of
$w=x^ay^b$ contains every element $z \ne -id$ with $\tr(z)=-2$.

\medskip
{\bf Case 6:} $\frac{q+1}{2} \nmid a$ and $\frac{q+1}{2} \nmid b$.

Let $x,y$ be some elements of order $q+1$ and note that $x^a \ne \pm
id$ and $y^b \ne \pm id$. Similarly to Case 5, observe that unless
$\tr(x^a)=\tr(y^b)=0$ for all elements $x,y$ of order $q+1$, one can
find two matrices $x,y$ of order $q+1$ satisfying $\tr(x^a) \ne
-\tr(y^b)$ (see \lemref{lem.traces}). Let $s=\tr(x)$ and $t=\tr(y)$.
According to the table in \propref{prop.summ.tr.surj},
$f_{a,b}(s,t) \ne 0$, and so \eqref{eq.tr.min.2} holds.

Now, the only case left is $q \equiv 3 \bmod 4$ and
$a=b=\frac{p(q^2-1)}{8}$ (see \remarkref{rem.obs.3.4}$(iv)$). In
\propref{prop.obs.xaya} we will show that in this case the image of
$w=x^ay^b$ contains no element $z \ne -id$ with $\tr(z)=-2$,
yielding the obstruction given in the proposition.
\end{proof}

\subsection{Obstructions for surjectivity of the word $w(x,y)=x^ay^b$}

\begin{remark}\label{rem.obs.3.4}
In the course of the proof of Propositions \ref{prop.ab2} and
\ref{prop.abmin2}, we need to consider the following special cases:

\begin{enumerate}\renewcommand{\theenumi}{\it \roman{enumi}}
\item $q$ is even, $|x^a|,|y^b| \in \{1,2,3\}$ for all $x,y \in
\SL(2,q)$, and $|x^a|=|y^b|=3$ if $|x|=|y|=q-1$:

Namely, $q=2^e$, $e$ is even, $\frac{q-1}{3}| a, \  (q+1)| a, \ 2|a$
(and similarly for $b$). Hence, $a$ and $b$ are multiples of
$$\lcm\left(\frac{q-1}{3},q+1,2\right) =  \frac{2(q^2-1)}{3},$$
namely, $a,b \in \{\frac{2(q^2-1)}{3}, \frac{4(q^2-1)}{3} \}$. By
\remarkref{rem.a.b.mod}, it is enough to consider the case $a = b =
\frac{2(q^2-1)}{3}.$

\item $q$ is even, $|x^a|,|y^b| \in \{1,2,3\}$ for all $x,y \in
\SL(2,q)$, and $|x^a|=|y^b|=3$ if $|x|=|y|=q+1$:

Namely, $q=2^e$, $e$ is odd, $\frac{q+1}{3}| a, \  (q-1)| a, \ 2|a$
(and similarly for $b$). Hence, $a$ and $b$ are multiples of
$$\lcm\left(\frac{q+1}{3},q-1,2\right) =  \frac{2(q^2-1)}{3},$$
namely, $a,b \in \{\frac{2(q^2-1)}{3}, \frac{4(q^2-1)}{3} \}$. By
\remarkref{rem.a.b.mod}, it is enough to consider the case $a = b =
\frac{2(q^2-1)}{3}.$

\item $q$ ia odd, $|x^a|,|y^b| \in \{1,2,4\}$ for all $x,y \in
\SL(2,q)$, and $|x^a|=|y^b|=4$ if $|x|=|y|=q-1$:

Namely, $q \equiv 1 \bmod 4$, $\frac{q-1}{4}| a, \ \frac{q+1}{2}| a,
\ p|a$ (and similarly for $b$). Hence,
$$a=b=\lcm\left(\frac{q-1}{4},\frac{q+1}{2},p\right) =
\frac{p(q^2-1)}{8}.$$

\item $q$ is odd, $|x^a|,|y^b| \in \{1,2,4\}$ for all $x,y \in
\SL(2,q)$, and $|x^a|=|y^b|=4$ if $|x|=|y|=q+1$:

Namely, $q \equiv 3 \bmod 4$, $\frac{q+1}{4}| a, \  \frac{q-1}{2}|
a, \ p|a$ (and similarly for $b$). Hence,
$$a=b=\lcm\left(\frac{q+1}{4},\frac{q-1}{2},p\right) =
\frac{p(q^2-1)}{8}.$$

\item $q$ is odd, $|x^a|,|y^b| \in \{1,2,3\}$ for all $x,y \in
\SL(2,q)$, and $|x^a|=|y^b|=3$ if $|x|=|y|=q-1$:

Namely, $q \equiv 1 \bmod 6$, $\frac{q-1}{3}| a, \  \frac{q+1}{2}|
a, \ p|a$ (and similarly for $b$). Hence $a$ and $b$ are multiples
of
$$\lcm\left(\frac{q-1}{3},\frac{q+1}{2},p\right)=
\begin{cases}
\frac{p(q^2-1)}{6} & \text{if } q \equiv 1 \bmod{12}\\
\frac{p(q^2-1)}{12} & \text{if } q \equiv 7 \bmod{12}
\end{cases} .$$

\item $q$ is odd, $|x^a|,|y^b| \in \{1,2,3\}$ for all $x,y \in
\SL(2,q)$, and $|x^a|=|y^b|=3$ if $|x|=|y|=q+1$:

Namely, $q \equiv 5 \bmod 6$, $\frac{q+1}{3}| a, \ \frac{q-1}{2}| a,
\ p|a$ (and similarly for $b$). Hence $a$ and $b$ are multiples of
$$\lcm\left(\frac{q+1}{3},\frac{q-1}{2},p\right)=
\begin{cases}
\frac{p(q^2-1)}{6} & \text{if } q \equiv 11 \bmod{12}\\
\frac{p(q^2-1)}{12} & \text{if } q \equiv 5 \bmod{12}
\end{cases} .$$
\end{enumerate}
\end{remark}


In order to investigate these obstructions in detail, we need the
following technical result on unipotent elements.

It follows from \thmref{thm.Mac} and \secref{sect.prop.SL} that for
any matrix $z \in \SL(2,q)$ with $\tr(z) \ne \pm 2$ and any two
integers $m,n>2$ dividing $p(q^2-1)$, one can find two matrices $x$
and $y$, such that $x^m = id = y^n$ and $z=xy$. However, a similar
result fails to hold if $z$ is \emph{unipotent}.

\begin{prop}\label{prop.prod.unip}
Let $z \in \SL(2,q)$ be a unipotent element (i.e. $z \ne \pm id$ and
$\tr(z) = \pm 2$), and let $m,n>2$ be two integers dividing
$p(q^2-1)$. Then there exist $x, y \in \SL(2,q)$, such that $x^m =
id = y^n$ and $z=xy$, if and only if none of the following
conditions hold:

\begin{enumerate}\renewcommand{\theenumi}{\it \roman{enumi}}
\item $q=2^e$, $e$ is odd, $\tr(z)=0$ and $m=n=3$;
\item $q \equiv 3 \bmod 4$, $\tr(z)=\pm 2$ and $m=n=4$;
\item $q \equiv 5 \bmod 6$, $\tr(z)=2$ and $m=n=3$.
\end{enumerate}
\end{prop}

\begin{proof}
If $m,n>2$ are two integers dividing $p(q^2-1)$, then one can find
$m',n'>2$ satisfying $m'|m$, $n'|n$ and moreover, either $m'=p$ or
$m'|q-1$ or $m'|q+1$, and either $n'=p$ or $n'|q-1$ or $n'|q+1$.
Thus, there exist some matrices $x,y$ in $\SL(2,q)$ such that $x$
has order $m'$ and $y$ has order $n'$, namely $x^{m'} = id =
x^{n'}$, and so $x^m = id = x^n$.

Assume that $\tr(z)=2$. If $m'=p$ then we can take $x=z$ and $y=id$.
Thus we may assume that both $m'$ and $n'$ are relatively prime to
$p$. Hence, unless $m'=n'=3$ or $m'=n'=4$, one can find two matrices
$x_1,y_1 \in \SL(2,q)$ such that $x_1^{m'} = id = y_1^{n'}$ and $\tr(x_1)
\ne \tr(y_1)$ (see \lemref{lem.traces}). Let $s=\tr(x_1)$ and
$t=\tr(y_1)$.

By \thmref{thm.Mac}, there exist two matrices $x_2,y_2 \in \SL(2,q)$
with $s=\tr(x_2), \ t=\tr(y_2)$ and $\tr(x_2y_2)=2$. Since $s \ne t$
then $x_2y_2 \ne id$, and hence, by \corref{cor.conj.tr}, there
exist some $x,y \in \SL(2,q)$ with $\tr(x)=s$ and $\tr(y)=t$
satisfying $z=xy$.

Now, assume that $\tr(z)=-2$. Unless $m'=n'=4$, one can find two
matrices $x_1,y_1 \in \SL(2,q)$ such that $x_1^{m'} = id = y_1^{n'}$ and
$\tr(x_1) \ne -\tr(y_1)$ (see \lemref{lem.traces}). Let $s=\tr(x_1)$
and $t=\tr(y_1)$. By \thmref{thm.Mac}, there exist two matrices
$x_2,y_2 \in \SL(2,q)$ with $s=\tr(x_2), \ t=\tr(y_2)$ and
$\tr(x_2y_2)=-2$. Since $s \ne -t$ then $x_2y_2 \ne -id$, and hence,
by \corref{cor.conj.tr}, there exist some $x,y \in \SL(2,q)$ with
$\tr(x)=s$ and $\tr(y)=t$ satisfying $z=xy$.

It is left to consider the cases $m'=n'=4$ and $m'=n'=3$.

For an odd $q$, in case $m'=n'=4$ we have $s=\tr(x)=\tr(y)=t=0$ and
$\om_t^2=-4.$ In \lemref{lem.yt} it is shown that such pair with
$x\ne y^{-1}$ exists if and only if $\om_t\in \BF_{q},$ therefore if
and only if  $q \equiv 1 \bmod 4.$

In case $m'=n'=3$ we have $s=\tr(x)=\tr(y)=t=-1$ and
$\om_t^2=-3.$ Hence, $\om_t\in \BF_{p^e}$ if and only if either $e$ is
even or $e$ is odd and $p \equiv 1 \bmod 6$, namely if and only if
$q=p^e \equiv 1 \bmod 6$.

In case $q=2^e$ and $m'=n'=3$ we have $s=\tr(x)=\tr(y)=t=1$ and
$\nu_{1,2}$ are the roots of the polynomial $\al^2+\al+1.$ These
roots belong to the field $\BF_{2^e}$ if and only if $e$ is even.
\end{proof}

The following proposition shows that the condition that neither
$a$ nor $b$ is divisible by the exponent of $\PSL(2,q)$ is not
sufficient for the surjectivity of the word $x^ay^b$ on $\PSL(2,q)$
(and on $\SL(2,q) \setminus \{-id\}$), yielding the obstructions given in
Propositions~\ref{prop.ab2} and~\ref{prop.abmin2}.

\begin{prop}\label{prop.obs.xaya}
Let $q$ be a prime power, $a,b \geq 1$, and $z \in \SL(2,q)$ a
unipotent element, satisfying the conditions given in the following
table.
\begin{center}
\begin{tabular}{c|c|c|c|c|}
  & $q=p^e$ & $a,b$ & conditions for $z$ & $m$\\
\hline \hline
$(i)$ & $q=2^e$, $e$ is odd & $a=b=\frac{2(q^2-1)}{3}$ & $z \ne id$ with $\tr(z)=0$ & $3$\\
\hline
$(ii)$ & $q \equiv 3 \bmod 4$ & $a=b=\frac{p(q^2-1)}{8}$ & $z \ne \pm id$ with $\tr(z) = \pm 2$ & $4$\\
\hline
$(iii)$ & $q \equiv 5 \bmod 12$ & $a,b \in \{\frac{p(q^2-1)}{6},\frac{p(q^2-1)}{12}\}$ & $z \ne id$ with $\tr(z) = 2$ & $3$\\
\hline
$(iv)$ & $q \equiv 11 \bmod 12$ & $a=b=\frac{p(q^2-1)}{6}$ & $z \ne id$ with $\tr(z) = 2$ & $3$\\
\hline
\end{tabular}
\end{center}
Then, in all these cases, $z$ does not belong to the image of $w=x^ay^b$.
\end{prop}

\begin{proof}
By \remarkref{rem.obs.3.4}, for every $x,y \in \SL(2,q)$ either
$x^a= \pm id$ or $x^a$ is of order $m$, and similarly for $y^b$. If
$z=x^ay^b$ is unipotent then necessarily $x^a\ne \pm id$ and $y^b
\ne \pm id$, hence both $x^a$ and $y^b$ are of order $m$. Assume
that $z$ is given as above. According to \propref{prop.prod.unip},
in all these cases $z$ cannot be written as a product of two
matrices of order $m$, hence, $z$ is not in the image of the word
map $w=x^ay^b$.
\end{proof}

\begin{prop}\label{prop.good.xaya}
Let $q$ be a prime power, $a,b \geq 1$, and $z \in \SL(2,q)$ a
unipotent element, satisfying the conditions given in the following
table.
\begin{center}
\begin{tabular}{c|c|c|c|c|}
& $q=p^e$ & $a,b$ & conditions for $z$ & $m$\\
\hline \hline
$(i)$ & $q=2^e$, $e$ is even & $a=b=\frac{2(q^2-1)}{3}$ & $z \ne id$ with $\tr(z)=0$ & $3$\\
\hline
$(ii)$ & $q \equiv 1 \bmod 4$ & $a=b=\frac{p(q^2-1)}{8}$ & $z \ne \pm id$ with $\tr(z) = \pm 2$ & $4$\\
\hline
$(iii)$ & $q \equiv 1 \bmod 12$ & $a=b=\frac{p(q^2-1)}{6}$ & $z \ne id$ with $\tr(z) = 2$ & $3$\\
\hline
$(iv)$ & $q \equiv 7 \bmod 12$ & $a,b\in \{\frac{p(q^2-1)}{6},\frac{p(q^2-1)}{12}\}$ & $z \ne id$ with $\tr(z) = 2$ & $3$\\
\hline
\end{tabular}
\end{center}
Then, in all these cases, $z$ is in the image of $w=x^ay^b$.
\end{prop}

\begin{proof}
According to  \propref{prop.prod.unip}  in all these cases there
exist two matrices of order $m$, $x_1$ and $y_1$, such that
$z=x_1y_1$. Moreover, by \remarkref{rem.obs.3.4}, any element $x$
order $q-1$ satisfies that $x^a$ has order $m$. Hence, there exists
some $x \in \SL(2,q)$ of order $q-1$ such that $x^a=x_1$ (see
\secref{sect.prop.SL}). Similarly, one can find some $y \in
\SL(2,q)$ of order $q-1$ such that $y^b = y_1$, and then $x^ay^b =
z$ as needed.
\end{proof}

\subsection{Missing $-id$ in the word map}\label{sect.min.id}


\begin{proof}[Proof of \thmref{thm.min.id}]
Assume that $q$ is odd and let $K = \max \left\{k: 2^k \ divides \ \frac{q^2-1}{2} \right\}.$

Observe that since $2^K \mid \frac{q^2-1}{2}$ and $\gcd(q-1,q+1)=2$,
then exactly one of the following holds:
\begin{itemize}
\item either $q-1 = 2^K \cdot m$ and $q+1 = 2 \cdot l$ for some odd
integers $l,m$;
\item or $q+1 = 2^K \cdot m$ and $q-1 = 2 \cdot l$ for some odd
integers $l,m$.
\end{itemize}

If $2^K \nmid a$ then one can write $a = 2^k a'$ for some $0 \leq
k<K$ and some odd integer $a'$. Without loss of generality we may
assume that $q-1 = 2^K \cdot m$ for some odd integer $m$.

Let $x_1 \in \SL(2,q)$ be some element of order $q-1$ and let $x=
x_1^{2^{K-k-1}}$, then
\[
    x^a = (x^{2^k})^{a'} = (x_1^{2^{K-1}})^{a'} = (x_1^{\frac{q-1}{2}})^{m
    a'} = (-id)^{ma'} = -id,
\]
and hence $-id = x^a id^b$ as needed.

On the other direction, if $2^K \mid a$ then since any element $x$
in $\SL(2,q)$ is either of order $p$ or of order dividing $q-1$ or
of order dividing $q+1$, we deduce that $x^a$ is either trivial or
of odd order. Similarly, if $2^K \mid b$ then for any element $y \in
\SL(2,q)$, $y^b$ is either trivial or of odd order.

If $-id=x^ay^b$ then neither $x^a$ nor $y^b$ is trivial. Let $l$ and
$m$ be the orders of $x^a$ and $y^b$ respectively, then both $l,m$
are odd and divide either $q-1$ or $q+1$. Without loss of generality
we may assume that both orders of $x$ and $y$ divide $q-1$, and that
$x^a$, and so also $y^b$, are in diagonal form, namely:
\[
x^a = \begin{pmatrix} \la & 0 \\ 0 & \la^{-1}\end{pmatrix}, \quad
y^b = \begin{pmatrix} \mu & 0 \\ 0 & \mu^{-1}\end{pmatrix},
\]
for some $\la,\mu \in \BF_q$ satisfying $\la^l=1$ and $\mu^m=1$.

Hence,
\[
-id = x^ay^b = \begin{pmatrix} \la\mu & 0 \\ 0 &
\la^{-1}\mu{-1}\end{pmatrix},
\]
implying that $\la\mu = -1$, but then, since $lm$ is odd,
\[
-1 = (-1)^{lm} = (\la\mu)^{lm} = (\la^l)^m(\mu^m)^l = 1\cdot 1 = 1,
\]
yielding a contradiction.
\end{proof}


\end{document}